\documentclass[12pt,reqno]{amsart}

\usepackage{amscd}
\usepackage{amssymb}
\usepackage{epsfig}
\usepackage{float}
\usepackage{url}
\newcommand{\Z}{{\mathbb Z}}
\newcommand{\Q}{{\mathbb Q}}
\newcommand{\C}{{\mathbb C}}
\newcommand{\R}{{\mathbb R}}
\renewcommand{\P}{{\mathbb P}}

\newcommand{\QQ}{{\mathcal Q}}

\newcommand{\www}{\widetilde}

\newcommand{\oooo}{\overline}
\newcommand{\uuuu}{\underline}

\newcommand{\paa}{\partial}

\DeclareMathOperator{\diag}{diag}

\DeclareMathOperator{\Eig}{Eig}

\DeclareMathOperator{\id}{id}
\DeclareMathOperator{\Imm}{Im}

\DeclareMathOperator{\rk}{rk}
\DeclareMathOperator{\Rad}{Rad}
\DeclareMathOperator{\Ree}{Re}
\DeclareMathOperator{\Seif}{Seif}
\DeclareMathOperator{\Sp}{Sp}
\DeclareMathOperator{\Spp}{Spp}

\begin{document}

\theoremstyle{plain}
\newtheorem{lemma}{Lemma}[section]
\newtheorem{definition/lemma}[lemma]{Definition/Lemma}
\newtheorem{theorem}[lemma]{Theorem}
\newtheorem{proposition}[lemma]{Proposition}
\newtheorem{corollary}[lemma]{Corollary}
\newtheorem{conjecture}[lemma]{Conjecture}
\newtheorem{conjectures}[lemma]{Conjectures}

\theoremstyle{definition}
\newtheorem{definition}[lemma]{Definition}
\newtheorem{withouttitle}[lemma]{}
\newtheorem{remark}[lemma]{Remark}
\newtheorem{remarks}[lemma]{Remarks}
\newtheorem{example}[lemma]{Example}
\newtheorem{examples}[lemma]{Examples}
\newtheorem{notations}[lemma]{Notations}
\newtheorem{recipe}[lemma]{Recipe}

\title[Spectral numbers for upper triangular matrices] 
{Conjectures on spectral numbers for upper triangular
matrices and for singularities}

\author{Sven Balnojan and Claus Hertling}

\address{Sven Balnojan\\
Lehrstuhl f\"ur Mathematik VI, Universit\"at Mannheim, 
Seminargeb\"aude A 5, 6, 68131 Mannheim, Germany}

\email{sbalnoja@mail.uni-mannheim.de}

\address{Claus Hertling\\
Lehrstuhl f\"ur Mathematik VI, Universit\"at Mannheim, 
Seminargeb\"aude A 5, 6, 68131 Mannheim, Germany}

\email{hertling@math.uni-mannheim.de}



\date{December 01, 2017}

\thanks{This work was supported by the DFG grant He2287/4-1
(SISYPH)}

\begin{abstract}
Cecotti and Vafa proposed in 1993 a beautiful idea how to
associate {\it spectral numbers} $\alpha_1,...,\alpha_n\in\R$ 
to real upper triangular $n\times n$ matrices $S$ with
1's on the diagonal and eigenvalues of $S^{-1}S^t$
in the unit sphere. Especially, 
$\exp(-2\pi i\alpha_j)$ shall be the eigenvalues of $S^{-1}S^t$. 

We tried to make their idea rigorous,
but we succeeded only partially. This paper fixes
our results and our conjectures. 
For certain subfamilies of matrices their idea
works marvellously, and there the spectral numbers fit well
to natural (split) polarized mixed Hodge structures.
We formulate precise conjectures saying how this should
extend to all matrices $S$ as above.

The idea might become relevant in the context of 
semiorthogonal decompositions in derived algebraic geometry.
Our main interest are the cases of Stokes like matrices
which are associated to holomorphic functions with isolated
singularities (Landau-Ginzburg models). 
Also there we formulate precise conjectures (which overlap
with expectations of Cecotti and Vafa).
In the case of the chain type singularities, we have
positive results.

We hope that this paper will be useful for further
studies of the idea of Cecotti and Vafa.
\end{abstract}

\maketitle

\tableofcontents

\setcounter{section}{0}

\section{Introduction, conjectures and results}\label{s1}
\setcounter{equation}{0}

\noindent
Cecotti and Vafa proposed in \cite{CV93} a beautiful idea 
how to associate to upper triangular matrices in 
\begin{eqnarray}\label{1.1}
T(n,\R)&:=&\{S=(s_{ij})\in M(n\times n,\R)\, |\, s_{ij}=0
\textup{ for }i>j,\\
&&\hspace*{1cm} s_{ii}=1,
S^{-1}S^t\textup{ has eigenvalues in }S^1\}
\nonumber
\end{eqnarray}
(with $n\in\Z_{\geq 1}$)
$n$ {\it spectral numbers} $\alpha_1,...,\alpha_n\in\R$
such that $e^{-2\pi i\alpha_1},...,e^{-2\pi i\alpha_n}$
are the eigenvalues of $S^{-1}S^t$.
Furthermore they claim to have an almost rigorous proof
that the recipe works and that in the case of 
Landau-Ginzburg models the spectral numbers of its 
Stokes matrices coincide with the true spectral numbers.

We consider the recipe as imcomplete and see serious gaps
in it and in the arguments that in the case of Landau Ginzburg
models the spectral numbers coincide.
We discuss this below. Still we find the idea fascinating.

This paper is the result of our efforts to make the recipe
work. We succeeded only partially. 
We have certain subspaces of $T(n,\R)$ where the recipe 
works and which are hopefully big enough to be useful
for an extension of the recipe to all of $T(n,\R)$.
Below we formulate precise conjectures and results. 
The recipe is as follows.

\begin{recipe}\label{t1.1}
Start with some matrix $S_1\in T(n,\R)$.
Choose a path from the unit matrix $E_n$ to $S_1$ within
$T(n,\R)$, i.e. a continuous map $S:[0,1]\to T(n,\R)$  
with $S(0)=E_n$ and $S(1)=S_1$. Now choose {\it in a natural
way} $n$ continuous functions
$\alpha_j:[0,1]\to\R$, $j\in\{1,...,n\}$, such that
$\alpha_j(0)=0$ and $e^{-2\pi i\alpha_1(r)},...,
e^{-2\pi i\alpha_n(r)}$ are the eigenvalues of 
$S(r)^{-1}S(r)^t$. Then 
$\alpha_1(1),...,\alpha_n(1)$ are defined to be the
spectral numbers of $S_1$. 
\end{recipe}

\begin{remarks}\label{t1.2}
(i) The recipe assumes that $T(n,\R)$ is connected.
Cecotti and Vafa conjecture this
\cite[first half of page 590]{CV93}, but have no proof for it.
Our conjecture \ref{t1.6} (a) below will imply this, 
but we also have no proof for it.
But even if $T(n,\R)$ is connected, the spectral numbers
might depend on the chosen path.

\medskip
(ii) Even if a path is given, 
it might happen that for some $r\in(0,1)$ several eigenvalues
of $S(r)^{-1}S(r)^t$ coincide. Then at this parameter $r$
one can exchange the continuations at $r$ of the functions $\alpha_j$ 
for these eigenvalues. 
Then in general it is unclear whether and how to make
a most natural choice and how to make the phrase
{\it in a natural way} in the recipe \ref{t1.1} precise.
This holds especially if $\alpha_i(r)-\alpha_j(r)\in 2\Z-\{0\}$.

\medskip
(iii) Cecotti and Vafa proposed in 
\cite[footnote 6 on page 583]{CV93} to choose the path such that
for $r\in(0,1)$ all eigenvalues of $S(r)^{-1}S(r)^t$ 
are different. This is within $T(n,\R)$ for most matrices not 
possible because the eigenvalue $-1$ has for all matrices
in $T(n,\R)$ even multiplicity because 
$\det(S(r)^{-1}S(r)^t)=1$.

\medskip
(iv) Only on the pages 589+590 in \cite{CV93}, it is demanded
that the path is within $T(n,\R)$, not yet on page 583.
But if one chooses a path which leaves $T(n,\R)$ there are
two problems.
The resulting spectral numbers might depend on the path.
And the arguments with $tt^*$-geometry for the coincidence
of the Stokes matrix spectral numbers with the
true spectral numbers of a Landau-Ginzburg model will not work
\cite[first half of page 590]{CV93}. Because of both problems
we restrict to the recipe with paths within $T(n,\R)$.
\end{remarks}

We have two subfamilies $T_{{\rm HOR}1}(n,\R)$ and 
$T_{{\rm HOR}2}(n,\R)\subset T(n,\R)$ for which the
recipe \ref{t1.1} works. The families will be presented
in section \ref{s4}, but here we give their crucial properties
and show how and why the recipe works for them.

\begin{theorem}\label{t1.3}
(a) The subspaces $T_{{\rm HOR}1}(n,\R)$ and 
$T_{{\rm HOR}2}(n,\R)\subset T(n,\R)$ which are defined
in definition \ref{t4.4} (a) satisfy the following properties.

\begin{list}{}{}
\item[($\gamma$)]
$T_{{\rm HOR}k}(n,\R)$ (for $k\in\{1,2\}$) 
can be represented by a closed simplex
(the convex hull of $\dim T_{{\rm HOR}k}(n,\R)+1$ many points)
in $\R^{\dim T_{{\rm HOR}k}(n,\R)}$. And 
\begin{eqnarray}\label{1.2}
\begin{array}{lll}
 & \dim T_{{\rm HOR}1}(n,\R) & \dim T_{{\rm HOR}2}(n,\R) \\
n\textup{ odd} & \frac{n-1}{2} & \frac{n-1}{2} \\
n\textup{ even} & \frac{n}{2} & \frac{n-2}{2} 
\end{array}
\end{eqnarray}
\item[($\beta$)]
For each $S\in T_{{\rm HOR}k}(n,\R)$, there is a regular matrix
$R^{mat}_{(k)}(S)\in GL(n,\R)$ with eigenvalues in $S^1$ and with
\begin{eqnarray}\label{1.3}
(-1)^k\cdot S^{-1}S^t &=& R^{mat}_{(k)}(S)^n.
\end{eqnarray}
{\rm Regular} means that $R^{mat}_{(k)}(S)$ has for each eigenvalue
only one Jordan block. The map $R^{mat}_{(k)}:T_{{\rm HOR}k}(n,\R)
\to GL(n,\R)$ is as a map to $M(n\times n,\R)$ affine linear.
\item[($\gamma$)]
$R^{mat}_{(k)}(S)$ is semisimple (and thus has pairwise different
eigenvalues) if and only if 
$S\in\textup{int}(T_{{\rm HOR}k}(n,\R))$.
\item[($\delta$)]
$E_n\in\textup{int}(T_{{\rm HOR}k}(n,\R))$ and
$R^{mat}_{(k)}(E_n)$ has the eigenvalues
$e^{-2\pi i(j-\frac{k}{2})/n}$, $j\in\{1,...,n\}$.
Furthermore, $\bigcap_{k=1,2}T_{{\rm HOR}k}(n,\R)=\{E_n\}$.

\end{list}

\medskip
(b) The recipe \ref{t1.1} works well within 
$T_{{\rm HOR}k}(n,\R)$. For 
$S_1\in T_{{\rm HOR}k}(n,\R)$ choose any continuous path
$S:[0,1]\to T_{{\rm HOR}k}(n,\R)$ with
$S(0)=E_n,S(1)=S_1$ and 
$S([0,1))\subset \textup{int}(T_{{\rm HOR}k}(n,\R))$.
Then for $r\in [0,1)$ the eigenvalues of 
$R^{mat}_{(k)}(S(r))$ are pairwise different and the paths 
$\alpha_1,...,\alpha_n:[0,1]\to T_{{\rm HOR}k}(n,\R)$
can be chosen uniquely such that $\alpha_j(0)=0$ and 
$e^{-2\pi i (\alpha_j(r)+j-\frac{k}{2})/n}$
for $j\in\{1,...,n\}$ are the eigenvalues of $R^{mat}_{(k)}(S(r))$.
The values $\alpha_1(1),...,\alpha_n(1)$ are independent
of the chosen path $S$ and give the spectrum
$\Sp(S)=\sum_{j=1}^n(\alpha_j(1))\in\Z_{\geq 0}(\R)$.
\end{theorem}

{\bf Proof:} Part (a) will be proved in section \ref{s4}.
Part (b) follows immediately from part (a).
In fact, part (a) implies existence and uniqueness
of continuous functions 
$\alpha_j^{(k)}:T_{{\rm HOR}k}(n,\R)\to \R$ such that
$\alpha_j^{(k)}(E_n)=0$ and 
$e^{-2\pi i(\alpha_j^{(k)}(S)+j-\frac{k}{2})/n}$
for $j\in\{1,...,n\}$ are the eigenvalues of 
$R^{mat}_{(k)}(S)$ for any $S\in T_{{\rm HOR}k}(n,\R)$.
For any $S\in T_{{\rm HOR}k}(n,\R)$ the values 
$\alpha_j^{(k)}(S)$ at $S$ are the spectral numbers of $S$.
The only matrix in $\bigcap_{k=1,2}T_{{\rm HOR}k}(n,\R)$
is $E_n$. Both cases $k=1$ and $k=2$ associate to $E_n$ 
the spectrum $\Sp(E_n)=\sum_{j=1}^n(0)$.
\hfill$\Box$

\begin{remarks}\label{t1.4}
(i) The crucial points are, that the matrices $R^{mat}_{(k)}(S)$ 
for $S\in\textup{int}(T_{{\rm HOR}k}(n,\R))$ have pairwise
different eigenvalues and the $\alpha_j^{(k)}(S)$ are determined
by these eigenvalues and that the values 
$e^{-2\pi i\alpha_j^{(k)}(S)}$ are the eigenvalues 
of $S^{-1}S^t$ because of \eqref{1.3}.

\medskip
(ii) HOR are the initials of the authors Horocholyn, Orlik and Randell of \cite{Ho17} and \cite{OR77}.
In \cite[ch. 2]{Ho17} half of the matrices in 
$\bigcup_{k=1,2}T_{{\rm HOR}k}(n,\R)$ were studied and the 
crucial equation \eqref{1.3} was proved for them. 
In \cite[(4.1) Conjecture]{OR77} it was conjectured that
special matrices $S$ in $\bigcup_{k=1,2}T_{{\rm HOR}k}(n,\Z)$
turn up as Stokes matrices of the chain type singularities
(sections \ref{s6} and \ref{s7}). 
The main result Theorem (2.11) in \cite{OR77} is that 
$(-1)^kS^{-1}S^t$ is a monodromy matrix for such a singularity.
\end{remarks}

That the recipe \ref{t1.1} works for the matrices in
$\bigcup_{k=1,2}T_{{\rm HOR}k}(n,\R)$ is good news.
It lead us to a number of conjectures and results
which form the contents of this paper. We hope that they
will be useful for a complete positive solution of
recipe \ref{t1.1}.

The rest of this introduction has two purposes.
It fixes some notions and proposes the conjectures
\ref{t1.6}, \ref{t1.7} and \ref{t1.9} which guide us
through all of the paper. And it explains the structure
of the paper and sketches some main results.

{\bf Section \ref{s4}} introduces the subfamilies
$T_{{\rm HOR}k}(n,\R)\subset T(n,\R)$ for $k\in\{1,2\}$
of HOR-matrices and proves theorem \ref{t1.3} (a).
And it adds more precise information, especially, that the
spectral pairs and the eigenspace decompositions of such 
a matrix give rise to a natural split polarized mixed 
Hodge structure.

{\bf Section \ref{s3}} prepares this. It introduces isomorphic
subspaces $T^{scal}_{{\rm HOR}k}(n,\R)$ and it formalizes
and studies the recipe
\begin{eqnarray}\label{1.4}
(\textup{eigenvalues of }R^{mat}_{(k)}(S)\mapsto
(\textup{spectral numbers }\alpha_1,...,\alpha_n
\textup{ of }S),
\end{eqnarray}
which is implicit in the proof of theorem \ref{t1.3} (b).
This is elementary, but worth to be studied for itself.
Properties of these spectral numbers give, combined
with conjecture \ref{t1.9} on the spectral numbers
of holomorphic functions, 
new features of these spectral numbers.
The recipe \eqref{1.4} will also be extended to a recipe
for spectral pairs $\Spp(S)=\sum_{j=1}^n(\alpha_j,k_j)
\in\Z_{\geq 0}(\R\times\Z)$. 

{\bf Section \ref{s2}} discusses spectral pairs from
an abstract point of view. This is elementary,
but must be provided. It also offers a review on 
the classification in \cite{BH17} of the
{\it Seifert form pairs} in definition \ref{t1.5}.
The notions in this definition are needed for the
conjectures \ref{t1.6} and \ref{t1.7}.

\begin{definition}\label{t1.5}
(a) A Seifert form pair $(H_\R,L)$ consists of a finite
dimensional real vector space $H_\R$ and a nondegenerate
bilinear form $L:H_\R\times H_\R\to \R$ (which is in general
neither symmetric or antisymmetric). Its {\it monodromy}
is the (unique) automorphism $M:H_\R\to H_\R$
with $L(Ma,b)=L(b,a)$ for $a,b\in H_\R$.

\medskip
(b) {\it Hermitian} Seifert form pairs are classified in 
\cite{Ne95}. The classification of {\it real} Seifert form 
pairs in \cite{BH17} is reviewed in section \ref{s2}.

\medskip
(c) Trivial lemma: Any matrix $S\in GL(n,\R)$ gives rise to the
Seifert form pair $\Seif(S):=(M(n\times 1,\R),L)$ with
$L(a,b):=a^t\cdot S^t\cdot b$. Its monodromy $M$
is given by $M(a)=S^{-1}S^t\cdot a$. 

\medskip
(d) We define the sets $\Seif(n),\Eig(n)$, the projection
$pr_{SE}$, and the maps $\Psi_{Seif}$ and $\Psi_{Eig}$
as follows.
\begin{eqnarray}
\Seif(n)&:=& \{\textup{isomorphism classes of Seifert form pairs }
(H_\R,L)\nonumber \\
&& \textup{with }\dim H_\R=n\textup{ and with eigenvalues
of the}\nonumber\\
&&\textup{monodromy }M\textup{ in }S^1\},\label{1.5}\\
\Eig(n)&:=& \{\textup{unordered tuples of numbers }
\lambda_1,...,\lambda_n\in S^1\}\nonumber \\
&:=& (S^1)^n/S_n,\label{1.6}
\end{eqnarray}
\begin{eqnarray}\label{1.7}
pr_{SE}:\Seif(n)\to \Eig(n),&& 
{}[(H_\R,L)]\mapsto (\textup{eigenvalues of }M),\hspace*{0.5cm}\\
\Psi_{Seif}:T(n,\R)\to  \Seif(n),&&
S\mapsto [\Seif(S)],\label{1.8}\\
\Psi_{Eig}:= pr_{SE}\circ \Psi_{Seif}&:&T(n,\R)\to \Eig(n).
\label{1.9}
\end{eqnarray}

(e) The group $G_{sign,n}:=\{\pm 1\}^n$ acts on $T(n,\R)$ by
conjugation,
\begin{eqnarray}\label{1.10}
(\varepsilon_1,...,\varepsilon_n):S\mapsto
\diag(\varepsilon_1,...,\varepsilon_n)\cdot S\cdot
\diag(\varepsilon_1,...,\varepsilon_n)
\end{eqnarray}
for $(\varepsilon_1,...,\varepsilon_n)\in G_{sign,n}$. 
The group $G_{sign,n}$ is called {\it sign group}.
Of course, the maps $\Psi_{Seif}$ and $\Psi_{Eig}$
are $G_{sign,n}$-invariant.

\medskip
(f) A {\it Seifert form stratum} in $T(n,\R)$ is a 
union of components of one fiber of $\Psi_{Seif}$ which
are permuted transitively by $G_{sign,n}$.
An {\it eigenvalue stratum} in $T(n,\R)$ is a union of
components of one fiber of $\Psi_{Eig}$ which are 
permuted transitively by $G_{sign,n}$. 
\end{definition}

\begin{conjecture}\label{t1.6}
(a) $T_{{\rm HOR}1}(n,\R)$ intersects each eigenvalue stratum
in $T(n,\R)$. 

\medskip
(b) If $S_1,S_2\in \bigcup_{k=1,2}T_{{\rm HOR}k}(n,\R)$
are in the same eigenvalue stratum of $T(n,\R)$ then 
$\Sp(S_1)=\Sp(S_2)$. 

\medskip
(c) If $S_1,S_2\in \bigcup_{k=1,2}T_{{\rm HOR}k}(n,\R)$
are in the same Seifert form stratum of $T(n,\R)$ then
$Spp(S_1)=\Spp(S_2)$.
\end{conjecture}

If it is true, conjecture \ref{t1.6} (a) implies that
$T(n,\R)$ is connected, 
conjecture \ref{t1.6} (a)+(b) gives spectral
numbers $\Sp(S)$ for any matrix $S\in T(n,\R)$,
and conjecture \ref{t1.6} (a)+(c) gives spectral pairs
for any matrix $S$ in a Seifert form stratum which is met
by $\bigcup_{k=1,2}T_{{\rm HOR}k}(n,\R)$.
But these are not all Seifert form strata, as remark \ref{t2.11}
(vii) and remark \ref{t5.3} (ii) will show.
Unfortunately, for the other Seifert form strata, 
we have no precise idea how to lift $\Sp(S)$ to $\Spp(S)$. 

\begin{conjecture}\label{t1.7}
Also for the matrices $S$ in the Seifert form strata which are
not met by $\bigcup_{k=1,2}T_{{\rm HOR}k}(n,\R)$,
$\Sp(S)$ lifts in a natural way to $\Spp(S)$.
\end{conjecture}

\begin{remarks}\label{t1.8}
(i) For odd $n$, $T_{{\rm HOR}1}(n,\R)$ and $T_{{\rm HOR}2}(n,\R)$
are mapped by suitable elements of the sign group $G_{sign,n}$
to one another. For odd $n$ conjecture \ref{t1.6} (a)
is equivalent to the analogous conjecture for 
$T_{{\rm HOR}2}(n,\R)$. But for even $n$ 
$\dim T_{{\rm HOR}2}(n,\R)=\dim T_{{\rm HOR}1}(n,\R)-1$,
and we expect that $T_{{\rm HOR}2}(n,\R)$ meets for large
enough $n$ some other Seifert form strata than
$T_{{\rm HOR}1}(n,\R)$.
\end{remarks}

In {\bf section \ref{s6}} we will review some facts
on holomorphic map germs
$f:(\C^{m+1},0)\to (\C,0)$ with an isolated singularity at 0
and on $M$-tame functions $f:X\to\C$ with $\dim X=m+1$.
Especially, we will discuss the following.
In both cases, there is a {\it Milnor number} $\mu=\mu(f)
\in\Z_{\geq 1}$. 
In both cases, there is a
$\textup{Br}_\mu\ltimes G_{sign,\mu}$ orbit of matrices 
$S\in T(\mu,\Z):=T(\mu,\R)\cap GL(\mu,\Z)$.
Here $\textup{Br}_\mu$ is the braid group with 
$\mu$ strings.
We call these matrices {\it Stokes matrices}.
Then $(-1)^{m+1}S^{-1}S^t$ is a matrix of the 
(classical global) monodromy. 
In both cases, there are $\mu$ spectral pairs
$\Spp(f)=\sum_{j=1}^\mu(\alpha_j(f),l(f))\in
\Z_{\geq 0}(\Q\times\Z)$
which come from natural mixed Hodge structures.
The first entries are the spectral numbers
$\Sp(f)=\sum_{j=1}^\mu(\alpha_j(f))\in
\Z_{\geq 0}(\Q)$.
In a suitable numbering, the spectral numbers satisfy the symmetry
$\alpha_j(f)+\alpha_{\mu+1-j}(f)=\frac{m-1}{2}$.

Building on the conjectures \ref{t1.6} and \ref{t1.7},
we have a conjecture which embraces the claim
of Cecotti and Vafa for Landau-Ginzburg models.

\begin{conjecture}\label{t1.9}
Suppose that the conjectures \ref{t1.6} and \ref{t1.7} are true.
Let $f$ be a holomorphic map germ
$f:(\C^{m+1},0)\to (\C,0)$ with an isolated singularity at 0
or an $M$-tame function $f:X\to\C$ with $\dim X=m+1$.
Then any Stokes matrix $S$ of $f$ satisfies
\begin{eqnarray}\label{1.11}
\Spp(S) = \Spp(f)-\Bigl(\frac{m-1}{2},m).
\end{eqnarray}
\end{conjecture}

The resulting equality $\Sp(S)=\Sp(S)-\frac{m-1}{2}$
is in the case of $M$-tame functions equivalent to the claim
in \cite{CV93} that in the case of the Landau-Ginzburg models
recipe \ref{t1.1} gives the central charges.
Though the equality of spectral pairs is slightly stronger.
And the case of an isolated hypersurface singularity is not
covered by Landau-Ginzburg models, except for the
quasihomogeneous singularities, they are $M$-tame on $\C^{m+1}$.

The results proved in the {\bf sections \ref{s5} and \ref{s7}}
can be summarized as follows.

\begin{theorem}\label{t1.10}
(a) (Section \ref{s5}) 
In the cases $n=2$ and $n=3$, the conjectures \ref{t1.6} and 
\ref{t1.7} and the conjecture \ref{t1.9}
for function germs are true.

\medskip
(b) (Section \ref{s7}) 
In the case of any chain type singularity $f(x_0,...,x_m)$, 
the matrix $S\in T_{{\rm HOR},k}(\mu,\Z)$ 
with $k\equiv m(2)$ which is considered in 
\cite[(4.1) Conjecture]{OR77}, satisfies $\Sp(S)=\Sp(f)-\frac{m-1}{2}$.
\end{theorem}

Theorem \ref{t1.10} (b) and conjecture (4.1) in \cite{OR77},
which says that the matrix $S$ there is a Stokes matrix
of $f$, imply conjecture \ref{t1.9} for the chain type
singularities. For them $\Spp(S)=\Spp(f)-(\frac{m-1}{2},m)$ and
$\Sp(S)=\Sp(f)-\frac{m-1}{2}$ are equivalent, as the monodromy is 
semisimple.

{\bf Section \ref{s8}} formulates some critic with
an explicit example
on some arguments in \cite{CV93} around recipe \ref{t1.1}
which use $tt^*$ geometry. And it offers some speculations
about approaches towards a positive solution of
recipe \ref{t1.1}.

We thank Duco van Straten and Martin Guest for discussions.

\section{Enhanced real Seifert form pairs and spectral pairs}
\label{s2}
\setcounter{equation}{0}

\noindent
Matrices in $T(n,\Z)$ turn up when one considers isolated hypersurface
singularities. There such a matrix encodes
the Seifert form on the Milnor lattice with respect to a distinguished basis.
See section \ref{s6}. Also spectral pairs turn up there. But their origin
is different, it is a natural polarized mixed Hodge structure with 
semisimple automorphism. 

Nemethi \cite{Ne95} studied in the singularity case
the relation between Seifert form and spectral pairs and proved
that the isomorphism class of the real Seifert form is (within the singularity
cases) equivalent to the spectral pairs modulo $2\Z\times\{0\}$.
We recover this below, see lemma \ref{t2.9}.

Here we will give a general abstract discussion of the
relationship between (abstract) Seifert forms and spectral pairs. 
This builds especially on \cite{BH17}. We will start with the classification
of real Seifert form pairs. We will define an enhancement of a real Seifert
form pair which includes spectral pairs. 
We will also discuss the triangular shape of the matrices
in $T(n,\R)$, it is related to a semiorthogonal decomposition.
This leads to a reformulation of the question how to associate spectral pairs
to matrices in $T(n,\R)$.

\begin{notations}\label{t2.1}
Throughout the whole paper, $H_K$ is a finite dimensional vector space
over a field $K$. If $H_\R$ is given, then $H_\C=H_\R\otimes_\R\C
=H_\R\oplus iH_\R$ is the complexification of $H_\R$.

If $L:H_K\times H_K\to K$ is a bilinear form then
two subspaces $V_1,V_2\subset H_K$ are {\it $L$-orthogonal}
if $L(V_1,V_2)=L(V_2,V_1)=0$.
And then the left and right orthogonal subspaces to a subspace
$U\subset H_K$ are $U^{L \perp}:=\{a\in H_K\, |\, L(a,U)=0\}$
and $U^{\perp R}:=\{b\in H_K\, |\, L(U,b)=0\}$. 

If $M:H_K\to H_K$ is an automorphism, then
$M_s,M_u,N:H_K\to H_K$ denote its semisimple,
its unipotent and its nilpotent part with
$M=M_sM_u=M_uM_s$ and $N=\log M_u,e^N=M_u$.
If $K=\C$, denote 
$H_\lambda:=\ker(M_s-\lambda\cdot \id):H_\C\to H_\C$,
$H_{\neq 1}:=\bigoplus_{\lambda\neq 1}H_\lambda$,
$H_{\neq -1}:=\bigoplus_{\lambda\neq -1}H_\lambda$.
\end{notations}

\begin{definition}\label{t2.2}
(a) A {\it Seifert form pair} is a pair $(H_\R,L)$ where
$L:H_\R\times H_\R\to\R$ is a nondegenerate bilinear form.
It is called {\it irreducible} if $H_\R$ does not split
into two nontrivial (i.e. both $\neq\{0\}$) 
$L$-orthogonal subspaces.

(b) The {\it monodromy} $M:H_\R\to H_\R$ of a Seifert form pair
$(H_\R,L)$ is the unique automorphism with
\begin{eqnarray}\label{2.1} 
L(Ma,b)=L(b,a) \qquad\textup{for all }a,b\in H_\R.
\end{eqnarray}
The {\it eigenvalues of a Seifert form pair} are the eigenvalues
of its monodromy.
Two bilinear forms $I_s$ and $I_a$ on $H_\R$ are defined by 
\begin{eqnarray}\label{2.2}
I_s(a,b)&:=& L(b,a)+L(a,b)=L((M+\id)a,b),   \\
I_a(a,b)&:=& L(b,a)-L(a,b)=L((M-\id)a,b),   \nonumber
\end{eqnarray}

\medskip
(c) An $S^1$-Seifert form pair is a Seifert form pair with
eigenvalues in $S^1$. 
\end{definition}

In \cite{BH17}, also four other bilinear forms $I_s^{(2)},I_a^{(2)},I_s^{(3)}$
and $I_a^{(3)}$ (on subspaces of $H_\R$) are associated to a Seifert form pair. 
Here we will use only $I_s$, and that only in the discussion of
the case $n=3$ in section \ref{s5}. 

Part (a) of the following theorem is an immediate consequence of
the calculation $L(Ma,Mb)=L(Mb,a)=L(a,b)$ and the definitions of
$I_s$ and $I_a$. The parts (b) and (c) give the classification of
Seifert form pairs and are proved in \cite[Theorem 2.5 and Theorem 2.9]{BH17}.

\begin{theorem}\label{t2.3}
(a) Let $(H_\R,L)$ be a Seifert form pair. 
The three bilinear forms $L$, $I_s$ and $I_a$ are monodromy invariant.
The radical of $I_s$ is $\ker(M+\id)$, the radical of $I_a$ is
$\ker(M-\id)$. 

\medskip
(b) Any Seifert form pairs splits into a direct and $L$-orthogonal sum
of irreducible Seifert form pairs. The splitting is unique up to
isomorphism.

\medskip
(c) The irreducible Seifert form pairs are given by the types
with the following names.
\begin{eqnarray}\label{2.3}
\Seif(\lambda,1,n,\varepsilon)&\textup{with}&
(\lambda=1\ \& \ n\equiv 1(2))\\
&\textup{or}& (\lambda=-1\ \&\ n\equiv 0(2)),\nonumber\\
\Seif(\lambda,2,n)&\textup{with}&
(\lambda=1\ \& \ n\equiv 0(2))\label{2.4}\\
&\textup{or}& (\lambda=-1\ \&\ n\equiv 1(2)),\nonumber\\
\Seif(\lambda,2,n,\zeta)
&\cong& \Seif(\oooo{\lambda},2,n,\oooo{\zeta})
\label{2.5}\\
&\textup{with}& \lambda,\zeta\in S^1-\{\pm 1\},
\zeta^2=\oooo{\lambda}\cdot (-1)^{n+1},\nonumber\\
\Seif(\lambda,2,n)&\textup{with}&
\lambda\in\R_{>1}\cup \R_{<-1},\label{2.6}\\
\Seif(\lambda,4,n)&\textup{with}&
\lambda\in\{\zeta\in \C\, |\, |\zeta|>1,\Imm\zeta>0\}.\label{2.7}
\end{eqnarray}
Here $n\in\Z_{\geq 1},\varepsilon\in\{\pm 1\}$.
The types are uniquely determined by the properties 
above of $\lambda$ and $n$ and the following properties.
$H_\lambda$ and Jordan blocks are meant with respect to $M$. 

\begin{list}{}{}
\item[\eqref{2.3}] $\Seif(\lambda,1,n,\varepsilon):$ 
$\dim H_\R=n$, $H_\C=H_\lambda$, 
one Jordan block, for each $a\in H_\R-\Imm N$
$$L(a,N^{n-1}a)\in \varepsilon\cdot\R_{>0}.$$
\item[\eqref{2.4}] $\Seif(\lambda,2,n):$
$\dim H_\R=2n$, $H_\C=H_\lambda$, 
two Jordan blocks of size $n$. 
\item[\eqref{2.5}] $\Seif(\lambda,2,n,\zeta):$
$\dim H_\R=2n$, $H_\C=H_\lambda\oplus H_{\oooo\lambda}$, 
two Jordan blocks of size $n$, for each $a\in H_\lambda-\Imm N$
$$L(a,N^{n-1}\oooo{a})\in \zeta\cdot\R_{>0}.$$
\item[\eqref{2.6}] $\Seif(\lambda,2,n):$
$\dim H_\R=2n$, $H_\C=H_\lambda\oplus H_{\lambda^{-1}}$, 
two Jordan blocks of size $n$. 
\item[\eqref{2.7}] $\Seif(\lambda,4,n):$
$\dim H_\R=4n$, $H_\C=H_\lambda\oplus H_{\lambda^{-1}}
\oplus H_{\oooo\lambda}\oplus H_{\oooo{\lambda}^{-1}}$, 
four Jordan blocks of size $n$.
\end{list}
\end{theorem}

In this paper, only $S^1$-Seifert form pairs will be relevant.
So we will not need the types in \eqref{2.6} and \eqref{2.7}.
Only for completeness sake, we have given the full classification.

The signature of $I_s$ will be useful in the case $n=3$ in section
\ref{s5} for determining the irreducible Seifert form pairs.

\begin{lemma}\label{t2.4}\cite[Lemma 2.10]{BH17}
The following table lists for the irreducible Seifert form pairs
in theorem \ref{t2.3} (c) the signature of $I_s$. 
\begin{eqnarray*}
\begin{array}{llc}
\textup{type of a Seifert} & \textup{form pair}
& \textup{signature of }I_s \\ \hline
\Seif(1,1,n,\varepsilon) 
&\textup{with }n\equiv\varepsilon(4)
& (\frac{n+1}{2},0,\frac{n-1}{2}) \\
\Seif(1,1,n,\varepsilon) 
&\textup{with }n\equiv-\varepsilon(4)  
& (\frac{n-1}{2},0,\frac{n+1}{2}) \\
\Seif(-1,1,n,\varepsilon) 
&\textup{with }n-1\equiv\varepsilon(4)
& \hspace*{0.3cm} (\frac{n}{2},1,\frac{n-2}{2}) \\
\Seif(-1,1,n,\varepsilon) 
&\textup{with }n-1\equiv-\varepsilon(4)
& (\frac{n-2}{2},1,\frac{n}{2}) \hspace*{0.4cm} \\
\Seif(1,2,n) 
&\textup{(with }n\equiv 0(2))
& (n,0,n) \\
\Seif(-1,2,n) &
\textup{(with }n\equiv 1(2))   
& (n-1,2,n-1) \\
\Seif(\lambda,2,n,\zeta\varepsilon) 
&\textup{with } n\equiv 0(2)
& (n,0,n)  \\ 
&  (\textup{and }\lambda\in S^1-\{\pm 1\}) &  \\
\Seif(\lambda,2,n,\zeta) 
&\textup{ with }  n\equiv 1(2) 
& (n-1,0,n+1)  \\
& (\textup{and }\lambda\in S^1-\{\pm 1\}) &  \\
\Seif(\lambda,2,n,-\zeta) 
&\textup{with }  n\equiv 1(2) 
& (n+1,0,n-1)  \\
& (\textup{and }\lambda\in S^1-\{\pm 1\}) &  \\
\Seif(\lambda,2,n) 
&\textup{with }
\lambda\in \R_{>1}\cup\R_{<-1} 
& (n,0,n)  \\
\Seif(\lambda,4,n) 
&\textup{with }  
\lambda\in\{\zeta\in\C| 
& (2n,0,2n)  \\
& |\zeta|>1,\Imm\zeta>0\} & 
\end{array}
\end{eqnarray*}
Here $n\in\Z_{\geq 1},\varepsilon\in\{\pm 1\}$, and 
in the lines 7--9 
$\zeta:=\frac{\oooo{\lambda}+1}{|\lambda+1|}\cdot i^{n+1}$.
\end{lemma}

Now we turn to spectral pairs, first in an elementary abstract setting.

\begin{definition}\label{t2.5}
(a) A {\it spectral pair} is a pair $(\alpha,k)\in\R\times\Z$.
An unordered tuple of $n$ spectral pairs is denoted by
\begin{eqnarray*}
\Spp=\sum_{(\alpha,k)\in\R\times \Z}d(\alpha,k)(\alpha,k)
\in\Z_{\geq 0}(\R\times\Z)\subset \Z(\R\times \Z)
\end{eqnarray*}
with $|\Spp|:=\sum_{(\alpha,k)}d(\alpha,k)=n$. Here $\Z(\R\times\Z)$
is the group ring over $\R\times \Z$. 
The number $d(\alpha,k)$ is the multiplicity of $(\alpha,k)$ as
a spectral pair. Any numbering of the $n$
spectral pairs gives $\Spp=\sum_{j=1}^n(\alpha_j,k_j)$.

\medskip
(b) (i) A {\it spectral pair ladder} (short: spp-ladder) 
consists of $l+1$ spectral pairs
\begin{eqnarray}\label{2.8}
(\alpha+k,m+l-2k)\quad\textup{with }k\in\{0,1,...,l\}.
\end{eqnarray}
Here $m\in\Z$ and $l\in\Z_{\geq 0}$. The numbers $m$ and $l$ are uniquely
determined by the spectral pair ladder. $l+1$ is its {\it length},
and $m$ is its {\it center}.
Its {\it first spectral pair} is the pair $(\alpha,m+l)$. 
Its {\it first spectral number} is $\alpha$. 
The spp-ladder is determined by $m$ and $l$ and its first spectral number.

(ii) The {\it partner spp-ladder} is the spp-ladder
\begin{eqnarray}\label{2.9}
(m-l-1-\alpha+k,m+l-2k)\quad\textup{with }k\in\{0,1,...,l\}.
\end{eqnarray}
It has the same length and center. 
The {\it distance} of an spp-ladder to its partner is $2\alpha+l+1-m$. 

(iii) A spp-ladder is {\it single} if it is its own partner, i.e. if 
the distance to its partner is 0, i.e. if $\alpha=\frac{m-l-1}{2}$. 

\medskip
(c) An unordered pair of spp-ladders (short: sppl-pair) consists
of two spp-ladders which are partners of one another and which
have distance $\neq 0$. 
\end{definition}

\begin{lemma}\label{t2.6}
(a) Each sppl-pair and each single spp-ladder 
with center $m$ are invariant under the
Kleinian group $\id,\pi_1,\pi_2,\pi_3:\R\times\Z\to\R\times\Z$ with
\begin{eqnarray}\label{2.10}
\pi_1:(\frac{m-1}{2}+\alpha,m+k)&\mapsto&
(\frac{m-1}{2}-\alpha,m-k),\\
\pi_2:(\frac{m-1-k}{2}+\alpha,m+k)&\mapsto&
(\frac{m-1-k}{2}-\alpha,m+k),\nonumber\\
\pi_3=\pi_1\circ\pi_2=\pi_2\circ\pi_1:
(\alpha,m+k)&\mapsto&
(\alpha+k,m-k).\nonumber
\end{eqnarray}
In the case of a sppl-pair, $\pi_3$ maps each spp-ladder to itself,
$\pi_1$ and $\pi_2$ map the two spp-ladders to one another.

(b) Suppose that a tuple $\Spp\in\Z_{\geq 0}(\R\times \Z)$ of $n$
spectral pairs is built from sppl-pairs and single spp-ladders with center $m$.
Then the sppl-pairs and the single spp-ladders are uniquely determined
by $\Spp$.
 
(c) Suppose that a tuple $\Spp\in\Z_{\geq 0}(\R\times \Z)$ of $n$
spectral pairs is built from sppl-pairs and single spp-ladders with center $m$.
Then $\Spp\mod 2\Z\times\{0\}$ determines each spp-ladder uniquely 
up to simultaneous shift of its members by elements of $2\Z\times\{0\}$,
so it determines the lengths, the centers and the first spectral numbers
modulo $2\Z$ of all spp-ladders.
\end{lemma}

{\bf Proof:} Trivial. \hfill$\Box$

\begin{definition}\label{t2.7}
Let $(H_\R,L)$ be an $S^1$-Seifert form pair.

(a) An {\it enhancement} of it is a decomposition 
of $(H_\R,L)$ into a direct and $L$-orthogonal sum of Seifert form
pairs $(H_\R^{(j)},L^{(j)})$ with $j\in\{1,...,r\}$ 
for some $r\in\Z_{\geq 1}$ together with
spectral pairs $\Spp^{(j)}\in\Z_{\geq 0}(\R\times\Z)$ with the 
following properties. 

\begin{list}{}{}
\item[(i)]
$\Spp^{(j)}$ consists of finitely many copies of the same sppl-pair 
or the same single spp-ladder. Its length is called $l_j$. 
All sppl-pairs and spp-ladders in $\Spp:=\sum_{j=1}^r\Spp^{(j)}$ 
have the same center $m\in\Z$.
This is also called the center of the enhancement.
The first spectral number of the/one of the two spp-ladders is called
$\alpha_j$ (if there are two, it does not matter which one).
$|\Spp^{(j)}|=\dim H_\R^{(j)}$.
\item[(ii)]
$(H_\R^{(j)},L^{(j)})$ decomposes into copies 
of one irreducible Seifert form pair 
\begin{eqnarray*}
\Seif((-1)^{m+1}e^{-2\pi i\alpha_j},2,l_j+1,\zeta_j)\textup{ in \eqref{2.5}} 
&\textup{if}& 2\alpha_j+l_j+1-m\in \R-\Z,\\
\Seif((-1)^{m+1}e^{-2\pi i\alpha_j},2,l_j+1)\textup{ in \eqref{2.4}}
&\textup{if}& 2\alpha_j+l_j+1-m\in \Z-2\Z,\\ 
\Seif((-1)^{m+1}e^{-2\pi i\alpha_j},1,l_j+1,\varepsilon_j)\textup{ in \eqref{2.3}}
&\textup{if}& 2\alpha_j+l_j+1-m\in 2\Z.
\end{eqnarray*}
\end{list}

(b) An enhancement with center $m$ is {\it polarized} if in (a)(ii)
\begin{eqnarray}\label{2.11}
(\varepsilon_j\textup{ resp. }\zeta_j)
=e^{\frac{1}{2}\pi i(2\alpha_j+l_j+1-m)}.
\end{eqnarray}
An enhancement with center $m$ is {\it signed polarized} if in (a)(ii)
\begin{eqnarray}\label{2.12}
(\varepsilon_j\textup{ resp. }\zeta_j)
=(-1)^{l_j} e^{\frac{1}{2}\pi i(2\alpha_j+l_j+1-m)}.
\end{eqnarray}
\end{definition}

\begin{remarks}\label{t2.8}
(i) Claim: An $S^1$-Seifert form pair $(H_\R,L)$ with (signed) polarized 
enhancement gives rise to  and is equivalent to a 
{\it split (signed) Steenbrink polarized mixed Hodge structure} on $H_\C$. 

The notions mixed Hodge structure and split mixed Hodge structure are
standard, see e.g. \cite[Def. 3.3 (a) and Remark 3.7]{BH17}.
The notion {\it Steenbrink polarized mixed Hodge structure}
is defined in \cite[Def. 3.3 (d)]{BH17}.
The signed version is defined in \cite[Def. 6.1]{BH17}.
The signed version turns up in the case of isolated hypersurface 
singualarities. The unsigned version turns up in $M$-tame functions.
For both cases see section \ref{s6}.

The claim follows easily from the results in \cite{BH17},
especially theorem 4.4. It builds on Deligne's $I^{p,q}$ of a mixed
Hodge structure, on the polarizing form of a polarized mixed Hodge structure,
and on the relation between Seifert form pairs and isometric triples,
which is developed in \cite{BH17}.

\medskip
(ii) Nemethi \cite{Ne95} considered the case of an isolated hypersurface 
singularity $f$ and studied there the relationship between the spectral pairs
$\Spp(f)$ of Steenbrink's mixed Hodge structure and the real Seifert form.
He found that $\Spp(f)\mod 2\Z\times\{0\}$ is equivalent to the isomorphism
class of the real Seifert form. The following lemma recovers this result
modulo the claim above in (i). 

\medskip
(iii) But for this result one has to know a priori that $\Spp(f)$ comes
from a signed Steenbrink polarized mixed Hodge structure,
or that $\Spp(f)$ is part of a {\it signed polarized} enhancement
of the real Seifert form.
\end{remarks}

\begin{lemma}\label{t2.9}
Two $S^1$-Seifert form pairs $(H_\R^{i},L^{i})$ for $i\in\{1,2\}$ 
with polarized enhancements (or with signed polarized enhancements)
with centers $m$ and spectral pairs $\Spp^i$ satisfy 
\begin{eqnarray}\label{2.13}
(H_\R^1,L^1)\cong (H_\R^2,L^2)\iff \Spp^1\equiv \Spp^2\mod 2\Z\times\{0\}.
\end{eqnarray}
\end{lemma}

{\bf Proof:} One can refine the decompositions of 
$(H_\R^1,L^1)$  and $(H_\R^2,L^2)$ 
in their enhancements to decompositions into sums of
irreducible Seifert form pairs such that each comes equipped with
a single spp-ladder or a sppl-pair.
Then the irreducible Seifert form pair determines 
the length $l$ of the single spp-ladder
or of each spp-ladder in the sppl-pair. The center of the
spp-ladder(s) is $m$.

The first spectral number $\alpha$ of the single spp-ladder 
or the first spectral numbers $\alpha$ and $\www\alpha$ of the 
two spp-ladders in the sppl-pair are determined modulo $\Z$ 
by $e^{-2\pi i\alpha}=(-1)^{m+1}\lambda$ and 
$e^{-2\pi i\www\alpha}=(-1)^{m+1}\oooo\lambda$ where $\lambda$
and $\oooo\lambda$ are the eigenvalue(s) of the irreducible Seifert form pair.

$\alpha$ and $\www\alpha$ are determined modulo $2\Z$ by the condition 
\eqref{2.11} respectively \eqref{2.12} in the cases \eqref{2.3} and \eqref{2.5}.
In the case \eqref{2.4}, they satisfy $\alpha\in\frac{1}{2}\Z$ and 
$\www\alpha\equiv\alpha+1(2)$.

Therefore the isomorphism class of $(H_\R^i,L^i)$ determines
the union $\Spp^i$ of all spp-ladders in the enhancement 
modulo $2\Z\times\{0\}$. This proves $\Rightarrow$. 

$\Leftarrow:$ Let $(\alpha_j,m_j,l_j)$ for $j\in\{1,...,\rho^1\}$ 
be the first spectral numbers, the centers and the lengths minus one 
of the spectral pair ladders in $\Spp^1$. By lemma \ref{t2.6} (c), 
the triples $(\alpha_j\mod 2\Z,m_j,l_j)$ are determined by 
$\Spp^1\mod 2\Z\times\{0\}$. Definition \ref{t2.7} 
and \eqref{2.3} and \eqref{2.4} show that each such triple 
determines a unique irreducible Seifert form pair
in $(H_\R^1,L^1)$. In the case of a sppl-pair, the triples of the two 
spp-ladders determine the same irreducible Seifert form pair.
This shows $\Leftarrow$. \hfill$\Box$

\bigskip
Finally, we put the matrices in $T(n,\R)$ into the frame
of Seifert form pairs.

\begin{lemma}\label{t2.10}
Let $(H_\R,L)$ be an $S^1$-Seifert form pair with $\dim H_\R=n\in\Z_{\geq 1}$. 
The following data are equivalent.
\begin{list}{}{}
\item[(A)]
A basis $\uuuu{v}=(v_1,...,v_n)$ with $L(\uuuu{v}^t,\uuuu{v})\in T(n,\R)$
up to the signs of the basis vectors $v_j$.
\item[(B)]
A splitting $H_\R=\bigoplus_{j=1}^n H_\R^{(j)}$ with
$\dim H_\R^{(j)}=1$, 
$L(H_\R^{(i)},H_\R^{(j)})=0$ for $i<j$ and 
$L(H_\R^{(j)},H_\R^{(j)})=\R_{\geq 0}$.
\item[(C)]
A complete flag $\{0\}\subset U_0\subset U_1\subset U_2
\subset ...\subset U_n=H_\R$
(complete flag means $\dim U_j=j$) with 
\begin{eqnarray}\label{2.14}
H_\R &=&\bigoplus_{j=1}^n H_\R^{(j)}\quad\textup{where }
H_\R^{(j)}:=U_j\cap U_{j-1}^{\perp R},\\
L(H_\R^{(j)},H_\R^{(j)})&=&\R_{\geq 0}.\label{2.15}
\end{eqnarray}
\end{list}
\end{lemma}

{\bf Proof:}
(A)$\Rightarrow$(B): Put $H_\R^{(j)}:=\R\cdot v_j$.

(B)$\Rightarrow$(A): For each $j$ choose a basis vector $v_j$ of
$H_\R^{(j)}$ with $L(v_j,v_j)=1$. It exists and is unique up to the sign.

(B)$\Rightarrow$(C): Put $U_j:=\bigoplus_{i\leq j}H_\R^{(i)}$.
Then $U_{j-1}^{\perp R}=\bigoplus_{i\geq j}H_\R^{(i)}$. 

(C)$\Rightarrow$(B): $H_\R^{(j)}$ has because of 
$\dim U_j+\dim U_{j-1}^{\perp R}=n+1$
at least dimension 1. By \eqref{2.14} it has dimension 1.\hfill$\Box$

\begin{remarks}\label{t2.11}
(i) A splitting as in (B) can be called a {\it semiorthogonal decomposition}.
Such splittings are considered in a much richer context in derived
algebraic geometry.

(ii) The complete flag in (C) and the positivity condition \eqref{2.15}
might remind one of Hodge structures. But there is no close relationship.

(iii) In the case of isolated hypersurfaces the data in lemma \ref{t2.10}
come from a distinguished basis, a refinement of the $\Z$-lattice structure. 
Steenbrink's mixed Hodge structure is of a transcendent origin and
has a clear relationship with the real structure, but no known
relationship with distinguished bases.

(iv) Nevertheless, the wish to associate to matrices $S\in T(n,\R)$ 
spectral pairs, can now be interpreted as the wish to see in the
data in lemma \ref{t2.10} a shadow of mixed Hodge structures.

(v) Let $(H_\R,L)$ be a real Seifert form pair. 
The set of all complete flags in  $H_\R$ is a real projective algebraic 
manifold $M^{flags}$.  
For any complete flag $U_\bullet$, the condition \eqref{2.14}
is equivalent to the condition
\begin{eqnarray}\label{2.16}
U_j\oplus U_{j}^{\perp R}=H_\R\quad\textup{for any }j\in\{1,...,n\}.
\end{eqnarray}
Let us call complete flags which do not satisfy \eqref{2.14} 
{\it degenerate}. They form a Zariski closed subvariety $M^{degen}$ in
$M^{flags}$, which separates the complement into components.
For each component a tuple $(\varepsilon_1,...,\varepsilon_n)\in\{\pm 1\}^n$
with
\begin{eqnarray}\label{2.17}
L(H_\R^{(j)},H_\R^{(j)})=\varepsilon_j\cdot \R_{\geq 0}
\end{eqnarray}
exists, where $U_\bullet$ is in the component and 
$H_\R^{(j)}$ is defined as in \eqref{2.14}.
This follows from the nondegeneracy of $L$. 
The components with $(\varepsilon_1,...,\varepsilon_n)=(1,...,1)$
give by (B)$\Rightarrow$(A) sets of matrices $L(\uuuu{v}^t,\uuuu{v})$ 
in $T(n,\R)$. The wish to associate to matrices $S\in T(n,\R)$
spectral pairs, is the wish to associate to each such component
spectral pairs.

(vi) A refinement of it is the wish to associate to each complete flag
in $M^{flags}-M^{degen}$ in a component with $(\varepsilon_1,...,\varepsilon_n)
=(1,...,1)$ an enhancement of $(H_\R,L)$. In the case of 
$S\in \bigcup_{k=1,2}T_{{\rm HOR}k}(n,\R)$, we will obtain such an enhancement.

(vii) There are Seifert form pairs $(H_\R,L)$ for which
$M^{flags}-M^{degen}$ has no components with 
$(\varepsilon_1,...,\varepsilon_n)=(1,...,1)$, 
i.e. which are not isomorphic to
$(M(n\times 1,\R),\www L)$ with $\www L(a,b)=a^t\cdot S^t\cdot b$
for any $S\in T(n,\R)$. Any sum of irreducible Seifert form pairs
\begin{eqnarray*}
\Seif(1,1,1,-1),\ \Seif(-1,1,2,-1),\ 
\Seif(-1,2,1),\ \Seif(\lambda,2,1,\zeta)
\end{eqnarray*}
(with $\lambda\in S^1-\{\pm 1\}$ and 
$\zeta=\frac{\oooo\lambda+1}{|\lambda+1|}\cdot i^{n+1}$)
has this property because then $I_s$ is negative (semi)definite
by lemma \ref{t2.4}.
In the cases $n\in\{2,3\}$ the only other Seifert form pairs
with this property are those which contain 
$\Seif(1,1,1,-1)$ or $\Seif(1,1,3,-1)$, 
see remark \ref{t5.3} (ii).
\end{remarks}

\section{A recipe for spectral pairs}
\label{s3}
\setcounter{equation}{0}

\noindent
Section \ref{s4} will present the subspaces 
$T_{{\rm HOR}k}(n,\R)$ of $T(n,\R)$ for $k\in\{1,2\}$ 
and study the properties of the matrices
in these subspaces. Here we prepare this.
We will introduce isomorphic subspaces
$T^{scal}_{{\rm HOR}k}(n,\R)\subset [0,1]^n\subset\R^n$
and $T^{pol}_{{\rm HOR}k}(n,\R)\subset \R[x]_{\deg=n}$ 
and propose for each of them a recipe for spectral pairs.

\begin{definition}\label{t3.1}
For $n\in\Z_{\geq 1}$ define the spaces 
\begin{eqnarray}\label{3.1}
T^{scal}_{{\rm HOR}1}(n,\R) &:=&
\{(\beta_1,...,\beta_n)\in[0,1]^n\, |\, 
\beta_1\leq ...\leq \beta_n,\\
&&\hspace*{2cm}\beta_j+\beta_{n+1-j}=1\},\nonumber \\
T^{scal}_{{\rm HOR}2}(n,\R) &:=&
\{(\beta_1,...,\beta_n)\in[0,1]^n\, |\, 
0=\beta_1\leq ...\leq \beta_n,\label{3.2}\\
&&\hspace*{2cm} \beta_j+\beta_{n+2-j}=1
\textup{ for }j\geq 2\},\nonumber \\
T^{simp}(n)&:=& \{(\beta_1,...,\beta_n)\in[0,\frac{1}{2}]^n
\, |\, \beta_1\leq ...\leq \beta_n\}. \label{3.3}
\end{eqnarray}
Define the map 
\begin{eqnarray}\label{3.4}
\Pi: \bigcup_{k=1,2}T^{scal}_{{\rm HOR}k}(n,\R)
&\to&\R[x]_{\deg=n}\\
\uuuu{\beta}=(\beta_1,...,\beta_n)&\mapsto& 
\prod_{j=1}^n (x-e^{-2\pi i\beta_j})\nonumber
\end{eqnarray}
and the spaces 
\begin{eqnarray}\label{3.5}
T^{pol}_{{\rm HOR}k}(n,\R):=\Pi(T^{scal}_{{\rm HOR}k}(n,\R))
\subset \R[x]_{\deg=n}\quad\textup{for }k\in\{1,2\}.
\end{eqnarray}
\end{definition}

\begin{lemma}\label{t3.2}
(a) $T^{simp}(n)$ is the $n$-simplex in $\R^n$ with the $n+1$
corners $(x_{1j},...,x_{nj})$ for $j\in\{0,1,...,n\}$ with
$x_{ij}=0$ for $i\leq j$ and $x_{ij}=\frac{1}{2}$ for 
$i>j$.

\medskip
(b) The following maps are affine linear isomorphisms. 
For odd $n$
\begin{eqnarray*}
T^{scal}_{{\rm HOR}1}(n,\R)\to T^{simp}(\frac{n-1}{2}),
&& (\beta_1,...,\beta_n)\mapsto 
(\beta_1,...,\beta_{\frac{n-1}{2}}),\\
T^{scal}_{{\rm HOR}2}(n,\R)\to T^{simp}(\frac{n-1}{2}),
&& (\beta_1,...,\beta_n)\mapsto 
(\beta_2,...,\beta_{\frac{n+1}{2}}).
\end{eqnarray*}
For even $n$ 
\begin{eqnarray*}
T^{scal}_{{\rm HOR}1}(n,\R)\to T^{simp}(\frac{n}{2}),
&& (\beta_1,...,\beta_n)\mapsto 
(\beta_1,...,\beta_{\frac{n}{2}}),\\
T^{scal}_{{\rm HOR}2}(n,\R)\to T^{simp}(\frac{n-2}{2}),
&& (\beta_1,...,\beta_n)\mapsto 
(\beta_2,...,\beta_{\frac{n}{2}}).
\end{eqnarray*}

(c) The map $\Pi$ in \eqref{3.4} is injective, and
\begin{eqnarray}
T^{pol}_{{\rm HOR}1}(n,\R)
&=& \{p\in \R[x]\, |\, \deg p=n,p_n=1, p_j=p_{n-j},\nonumber\\
&& \hspace*{2cm} 
\textup{all zeros of }p\textup{ are in }S^1\},\label{3.6}\\
T^{pol}_{{\rm HOR}2}(n,\R)
&=& \{p\in \R[x]\, |\, \deg p=n,p_n=1, p_j=-p_{n-j},\nonumber\\
&& \hspace*{2cm} 
\textup{all zeros of }p\textup{ are in }S^1\}.\label{3.7}
\end{eqnarray}
If $p\in T^{pol}_{{\rm HOR}k}(n,\R)$ then $p_0=(-1)^{k-1}$,
$p_j=p_0p_{n-j}$, $x^np(x^{-1})=p_0\cdot p(x)$,
$\lambda\in S^1$ and $\oooo{\lambda}$ have the same 
multiplicity as zeros of $p$, and the multiplicity
of 1 as a zero of $p$ is even for $k=1$ and odd for
$k=2$.
\end{lemma}

{\bf Proof:}
(a) Trivial. 

(b) For $\uuuu{\beta}\in T^{scal}_{{\rm HOR}1}(n,\R)$
the symmetry $\beta_j+\beta_{n+1-j}$ is used. 
For odd $n$ it implies $\beta_{\frac{n+1}{2}}=\frac{1}{2}$.
For $\uuuu{\beta}\in T^{scal}_{{\rm HOR}2}(n,\R)$
$\beta_1=0$ and the symmetry $\beta_j+\beta_{n+2-j}=1$
for $j\geq 2$ are used. For even $n$ the symmetry
implies $\beta_{\frac{n+2}{2}}=\frac{1}{2}$.

(c) Trivial. \hfill$\Box$ 

\bigskip
The following recipe formalizes the recipe 
\begin{eqnarray*}
\Bigl(\textup{eigenvalues of }R^{mat}_{(k)}(S)\Bigr)\mapsto 
\Bigl(\textup{spectral numbers }\Sp(S)\Bigr)
\end{eqnarray*}
which is implicit in the proof of theorem \ref{t1.3} (b).
In definition \ref{t4.4} (c) and theorem \ref{t4.5} (d) 
it is connected with theorem \ref{t1.3} (b).
It is completely elementary, but interesting in its own right.

\begin{recipe}\label{t3.3}
(a) The following recipe associates to any tuple
$\uuuu{\beta}\in T^{scal}_{{\rm HOR}k}(n,\R)$ for $k\in\{1,2\}$
a spectrum 
$\Sp(\uuuu\beta)=\sum_{j=1}^n(\alpha_j)\in\Z_{\geq 0}(\R)$.
Define for $j\in\{1,...,n\}$
\begin{eqnarray}\label{3.8}
\gamma_j&:=&\frac{1}{n}(j-\frac{k}{2})= \left\{\begin{array}{ll}
\frac{1}{n}(j-\frac{1}{2}) & \textup{ if }k=1,\\
\frac{1}{n}(j-1) & \textup{ if }k=2, \end{array}\right. \\
\alpha_j&:=& n(\beta_j-\gamma_j)=\left\{\begin{array}{ll}
n\beta_j-j+\frac{1}{2} & \textup{ if }k=1,\\
n\beta_j-j+1 & \textup{ if }k=2. \end{array}\right. \label{3.9}
\end{eqnarray}

(b) The following extends the recipe in (a) to a recipe for 
spectral pairs $\Spp(\uuuu\beta):=\sum_{j=1}^n(\alpha_j,k_j)
\in\Z_{\geq 0}(\R\times\Z)$. See lemma \ref{t3.4}
for the properties of $\Spp(\uuuu\beta)$. 
Consider $\kappa\in S^1$ with 
$\{\beta_j\, |\, e^{-2\pi i\beta_j}=\kappa\}\neq\emptyset$.
Then the recipe in (a) gives in fact
\begin{eqnarray}\label{3.10}
\sum_{j:\, \exp{(-2\pi i\beta_j)}=\kappa}
(\alpha_j)=\sum_{j=0}^l(\alpha+j)\textup{ for some }
\alpha\in\R,l\in\Z_{\geq 0}
\end{eqnarray}
(if $k=1$ and $\beta_1=0$ then 
$(\alpha_1,\alpha_n)=(\frac{-1}{2},\frac{1}{2})$, 
and if $k=2$ and $\beta_2=0$ then
$(\alpha_1,\alpha_2,\alpha_n)=(0,-1,1)$). 
Extend this to the spp-ladder $\sum_{j=0}^l(\alpha+j,1+l-2j)$
of length $l+1$ and center $m=1$ as in definition \ref{t2.5} (b),  
and define $\Spp(\uuuu\beta)$ as the sum of these spp-ladders.

\medskip
(c) For a polynomial $p\in T^{pol}_{{\rm HOR}k}(n,\R)$ define
the spectrum and the spectral pairs as follows,
\begin{eqnarray}\label{3.11}
\Sp(p):=\Sp(\Pi^{-1}(p)),\quad
\Spp(p):=\Spp(\Pi^{-1}(p)).
\end{eqnarray}
\end{recipe}

The spectral numbers $\alpha_1,...,\alpha_n$ in this recipe
are usually not ordered by size. But they satisfy the 
symmetry in part (b) of the following lemma.
The lemma states also properties of the spectral pairs.

\begin{lemma}\label{t3.4}
(a) Denote $\uuuu{\gamma}:=(\gamma_1,...,\gamma_n)$ in both
cases $k=1$ and $k=2$. Then
\begin{eqnarray}\label{3.12}
\uuuu{\gamma}\in T^{scal}_{{\rm HOR}k}(n,\R),\quad
\Pi(\uuuu{\gamma}) = x^n-(-1)^k,\\
\Spp(\uuuu\gamma)=n\cdot (0,1),\quad \Sp(\uuuu\gamma)
=n\cdot (0).\label{3.13}
\end{eqnarray}

(b) The spectral numbers $\alpha_1,...,\alpha_n$ in the 
recipe satisfy the symmetry
\begin{eqnarray}\label{3.14}
\alpha_j+\alpha_{n+1-j}&=& 0\quad\textup{for }k=1,\\
\alpha_1=0,\alpha_j+\alpha_{n+2-j}&=& 0\quad\textup{for }
k=2\textup{ and }j\geq 2.\label{3.15}
\end{eqnarray}
$\Spp(\uuuu\beta)$ consists of sppl-pairs and single spp-ladders
with center $m=1$, for each value $\kappa\in S^1$
with $\{\beta_j\ |\, e^{-2\pi i\beta_j}=\kappa\}\neq\emptyset$ 
one spp-ladder. The partner of the spp-ladder from $\kappa$
is the one from $\oooo\kappa$. 
The single spp-ladders are those which come from $\kappa\in\{\pm 1\}$, 
so there are at most two of them.

\medskip
(c) If $p\in T^{pol}_{{\rm HOR}k}(n,\R)$ then
$(-1)^np(-x)\in T^{pol}_{{\rm HOR}\www k}(n,\R)$ with
$\www k\equiv k+n(2)$, and then
\begin{eqnarray}\label{3.16}
\Spp(p)=\Spp((-1)^np(-x)).
\end{eqnarray}
\end{lemma}

{\bf Proof:}
(a) Trivial.

(b) $\uuuu\beta$ and $\uuuu\gamma$ are both in 
$T^{scal}_{{\rm HOR}k}(n,\R)$ and satisfy the same symmetry
in \eqref{3.1} or \eqref{3.2}.
Thus the tuple $\frac{1}{n}\uuuu\alpha=\uuuu\beta-\uuuu\gamma$
and the tuple $\uuuu\alpha$ satisfy the symmetry in 
\eqref{3.14} or \eqref{3.15}.

Consider as in part (b) of the recipe \ref{t3.3}
$\kappa\in S^1$ with $\{\beta_j\, |\, e^{-2\pi i\beta_j}=\kappa\}
\neq\emptyset$ and its spp-ladder. 
One sees easily with the symmetries \eqref{3.14} and \eqref{3.15}
that the spp-ladders for
$\kappa$ and $\oooo\kappa$ are partners.
Especially, those for $\kappa\in\{\pm 1\}$ are 
single spp-ladders.

(c) Write 
\begin{eqnarray*}
\www{p}(x):=(-1)^np(-x),\ \uuuu\beta:=\Pi^{-1}(p),
\ \uuuu{\www{\beta}}=\Pi^{-1}(\www p),
\ \www{p}\in T^{pol}_{{\rm HOR}\www{k}}(n,\R).
\end{eqnarray*}
Then $p_0=(-1)^{k-1}$ and $\www{p}_0=(-1)^{n+k-1}$ show
the first line of (c).

For even $n$
\begin{eqnarray*}
\www\beta_{\frac{n}{2}+j}=\frac{1}{2}+\beta_j
\quad\textup{ and }\quad\www\beta_j=-\frac{1}{2}+\beta_{\frac{n}{2}+j}
\textup{ for }j=1,...,\frac{n}{2}.
\end{eqnarray*}
For odd $n$ and $k=1$ 
\begin{eqnarray*}
\www\beta_{\frac{n+1}{2}+j}=\frac{1}{2}+\beta_j
\ \textup{ for }\ j=1,...,\frac{n-1}{2},\\
\www\beta_j=-\frac{1}{2}+\beta_{\frac{n-1}{2}+j}
\ \textup{ for }\ j=1,...,\frac{n+1}{2}.
\end{eqnarray*}
For odd $n$ and $k=2$ 
\begin{eqnarray*}
\www\beta_{\frac{n-1}{2}+j}=\frac{1}{2}+\beta_j
\ \textup{ for }\ j=1,...,\frac{n+1}{2},\\
\www\beta_j=-\frac{1}{2}+\beta_{\frac{n+1}{2}+j}
\ \textup{ for }\ j=1,...,\frac{n-1}{2}.
\end{eqnarray*}
Observe that 
\begin{eqnarray*}
\uuuu{\www\gamma} =_{def} \Pi^{-1}(\www{x^n-(-1)^k})
\stackrel{!}{=}\Pi^{-1}(x^n-(-1)^{\www k})
\end{eqnarray*}
is the $\uuuu\gamma$-vector for $\www k$. As $\uuuu{\www\gamma}$
is obtained from $\uuuu\gamma$ as any $\uuuu{\www\beta}$ 
from $\uuuu\beta$, the tuples of differences
$\uuuu{\www\beta}-\uuuu{\www\gamma}$ and 
$\uuuu\beta-\uuuu\gamma$ coincide up to reordering.
Therefore $\Sp(\www p)=\Sp(p)$. 
Its extension to $\Spp(\www p)=\Spp(p)$ is rather
obvious. \hfill$\Box$

\bigskip
It is interesting to ask about the images in 
$\Z_{\geq 0}(\R\times \Z)$ and in $\Z_{\geq 0}(\R)$
of the maps $\Spp$ and $\Sp$ from $T^{pol}_{{\rm HOR}k}(n,\R)$.
The answer is not difficult, it is given in the following
corollary. We omit the rather trivial proof.

\begin{corollary}\label{t3.5}
An unordered tuple $\sum_{\alpha\in\R}d(\alpha)(\alpha)
\in\Z_{\geq 0}(\R)$ of $n$ numbers 
(so $\sum_{\alpha\in\R}d(\alpha)=n$) is in 
$\Sp(T^{pol}_{{\rm HOR}k}(n,\R))$ if and only if the numbers
can be ordered as $\alpha_1,...,\alpha_n$ such that 
the symmetry in \eqref{3.14} respectively \eqref{3.15}
holds and $\alpha_{j+1}\geq \alpha_j-1$,
and in the case $k=1$ also $\alpha_1\geq -\frac{1}{2}$.

An unordered tuple $\sum_{(\alpha,k)\in\R\times\Z}
d(\alpha,k)(\alpha,k)\in\Z_{\geq 0}(\R\times \Z)$ of $n$ 
pairs is in $\Spp(T^{pol}_{{\rm HOR}k}(n,\R))$
if and only if the pairs can be ordered as
$(\alpha_1,k_1),...,(\alpha_n,k_n)$ such that
the conditions above hold and the tuple 
$\sum_{j=1}^n(\alpha_j,k_j)$ is obtained from the
tuple $\sum_{j=1}^n(\alpha_j)$ by part (b) of 
recipe \ref{t3.3}.
\end{corollary}

\begin{remarks}\label{t3.6}
(i) Corollary \ref{t3.5} implies that there is no gap
of size $>1$ in the spectral numbers if one orders them by size.
This follows (with the order in corollary \ref{t3.5}) 
from $\alpha_{j+1}\geq \alpha_j-1$, from $\alpha_1\geq -\frac{1}{2}$,
$\alpha_n\leq \frac{1}{2}$ for $k=1$, and from 
$\alpha_1=0$, $\alpha_2\geq -1$, $\alpha_n\leq 1$ for $k=2$.

\medskip
(ii) Now conjecture \ref{t1.9} implies that the spectral numbers of
isolated hypersurface singularities and $M$-tame functions
have no gap of size $>1$. 
This is not a very strong claim in the case of isolated
hypersurface singularities
(there usually the gaps between spectral numbers are much smaller),
but it is new in any case.
\end{remarks}

\begin{examples}\label{t3.7}
(i) It is also interesting to ask about the preimages
in $T^{pol}_{{\rm HOR}k}(n,\R)$ of spectral numbers or
spectral pairs, especially for $\Sp(f)$ with $f$ an
isolated hypersurface singularity or an M-tame function.
In these cases $\Spp(f)\in\Z_{\geq 0}(\Q\times \Z)$,
and, even stronger, the characteristic polynomial
$p_{ch,M}(x):=\prod_{j=1}^\mu(x-e^{-2\pi i \alpha_j})$ is in 
$\Z[x]$, i.e. it is a product of cyclotomic polynomials.

\medskip
(ii) In most cases, the preimages, the polynomials
$p\in T^{pol}_{{\rm HOR}k}(\mu,\R)$ with the correct
spectrum $Sp(p)=\Sp(f)-\frac{m-1}{2}$, are not in $\Z[x]$. 
If one looks only at the correct eigenvalues, and not at
the correct spectral numbers, one obtains the possibly
bigger set
\begin{eqnarray}\label{3.17}
\{p \in T^{pol}_{{\rm HOR}k}(\mu,\R)
&|& p(x)=\prod_{j=1}^\mu(x-\kappa_j),\\
&& p_{ch,M}(x)=\prod_{j=1}^\mu(x-(-1)^{k+m-1}\kappa_j^\mu)\}.
\nonumber
\end{eqnarray}
Even this set does often not contain polynomials in $\Z[x]$,
for example for the singularity $E_6$, see below (v).

\medskip
(iii) Remarkable exceptions are the chain type singularities,
which are treated in section \ref{s7}.
For them distinguished polynomials $p\in\Z[x]$ with
the correct spectrum $\Sp(p)=Sp(f)-\frac{m-1}{2}$ exist.
This will be proved in theorem \ref{t7.6}. The polynomials
$p$ are given in \eqref{7.6}.
In fact, in the moment, the chain type singularities
are the only candidates within isolated hypersurface 
singularities for which we know polynomials $p$ in 
$\Z[x]\cap T^{pol}_{{\rm HOR}k}(\mu,\R)$
with the correct spectrum.

\medskip
(iv) If $f(x_0,x_1)$ (so $m=1$) is one of the ADE-singularities, then the spectral numbers satisfy
$\frac{-1}{2}<\alpha_1\leq ...\leq \alpha_\mu <\frac{1}{2}$. 
Then the number of 
$\uuuu\beta\in T^{scal}_{{\rm HOR}k}(\mu,\R)$ with
$\Sp(\uuuu\beta)=\Sp(f)$ is (here $(2N)!!:=2^NN!$)
\begin{eqnarray*}
\mu!!&\textup{if}& \mu\textup{ is even and the singularity
 is not }D_\mu,\\
(\mu-1)!! &\textup{if}& \mu\textup{ is odd},\\
\mu!!\cdot\frac{1}{2}&\textup{if}&\textup{the singularity is }D_\mu 
\textup{ and }\mu\textup{ is even.}
\end{eqnarray*}
The numbers $\beta_j$ must satisfy
\begin{eqnarray}\label{3.18}
&&\beta_j = \gamma_j+\frac{1}{\mu}\alpha_{\sigma(j)},\\
&&\textup{the symmetry in }\eqref{3.1}\textup{ or }\eqref{3.2}, 
\textup{ including }\beta_1=0\textup{ in }\eqref{3.1},
\nonumber\\
&&0\leq \beta_1\leq ...\leq \beta_\mu\leq 1,\nonumber
\end{eqnarray}
here $\sigma\in S_\mu$ is a permutation.
Because of 
\begin{eqnarray*}
\max_j|\alpha_j|<\frac{1}{2}=1-\mu\gamma_\mu 
=\frac{\mu}{2}(\gamma_j-\gamma_{j-1})
=\mu \gamma_1-0,
\end{eqnarray*}
one can choose $\sigma\in S_\mu$ almost arbitrarily.
Only the symmetry in \eqref{3.1} or \eqref{3.2}
has to be observed.
For all ADE-singularities except $D_\mu$ with $\mu$ even,
the spectral numbers are pairwise different. 
For $D_\mu$ with $\mu$ even, 
$\alpha_{\frac{\mu}{2}}=\alpha_{\frac{\mu+2}{2}}=0$.

Though most of the polynomials $p=\Pi(\uuuu\beta)$ 
are not in $\Z[x]$. The singularities $A_\mu,D_\mu$ and $E_7$
can be written as chain type singularities.
Therefore by theorem \ref{t7.6} at least the polynomial in 
\eqref{7.6} is in $\Z[x]$. 

\medskip
(v) But the singularity $E_6$ and many other singularities
have characteristic polynomials $p_{ch,M}$ such that not even 
the set in \eqref{3.17} contains any polynomial
in $\Z[x]$. For $E_6$ as a curve singularity (so $m=1$) 
$p_{ch,M}=\Phi_{12}\Phi_6$. For $E_8$ as a curve singularity
$p_{ch,M}=\Phi_{15}$, and the set in \eqref{3.17} with $k=2$
contains the polynomial $p=\Phi_{15}\in \Z[x]$. But 
$\Sp(p)\neq \Sp(f)$. 
\end{examples}

\section{HOR-matrices}\label{s4}
\setcounter{equation}{0}

\noindent
In this section we will introduce the two subspaces
$T_{{\rm HOR}k}(n,\R)$ for $k\in\{1,2\}$ of $T(n,\R)$
and study the properties of matrices in these spaces.
We call the matrices HOR-matrices because of the 
initials of the authors Horocholyn, Orlik and Randell
of \cite{Ho17} and \cite{OR77}.
Horocholyn studied half of the matrices and proved
the crucial formula \eqref{4.20} \cite[ch. 2]{Ho17}.
Orlik and Randell considered a subfamily which is related
to the chain type singularities 
\cite[(4.1) Conjecture]{OR77} and which we will treat
in section \ref{s7}.
Before coming to the HOR-matrices, we recall a well
known fact from Picard-Lefschetz theory.
For the convenience of the reader, we present also a proof.

\begin{theorem}\label{t4.1}
Let $n\in\Z_{\geq 1}$, 
let $H_\R$ be an $\R$-vector space with a basis
$\uuuu{e}=(e_1,...,e_n)$, 
and let $S\in GL(n,\R)$.

(a) The matrix $S$ defines on $H_\R$ 
a bilinear form $L$, which is called {\rm Seifert form},
a symmetric bilinear form $I_s$,
an antisymmetric bilinear form $I_a$ and an automorphism
$M$, which is called {\rm monodromy}, by the 
formulas
\begin{eqnarray}\label{4.1}
L(\uuuu{e}^t,\uuuu{e})&=& S^t,\\
I_s(\uuuu{e}^t,\uuuu{e})&=& S+S^t,\quad
\textup{so }I_s(a,b)=L(b,a)+L(a,b),\label{4.2}\\
I_a(\uuuu{e}^t,\uuuu{e})&=& S-S^t,\quad
\textup{so }I_a(a,b)=L(b,a)-L(a,b),\label{4.3}\\  
M\, \uuuu{e}&=& \uuuu{e}\cdot S^{-1}S^t,\quad
\textup{so }L(Ma,b)=L(b,a).\label{4.4}
\end{eqnarray}
$L$ determines $I_s,I_a$ and $M$. 
The monodromy $M$ respects all three bilinear forms 
$L$, $I_s$ and $I_a$.

(b) Define endomorphisms $s^{(1)}_a$ and $s^{(2)}_b$
on $H_\R$ for $a\in H_\R$ with $I_s(a,a)=2$ and 
for arbitrary $b\in H_\R$ by
\begin{eqnarray}\label{4.5}
s^{(1)}_a(c):= c-I_s(a,c)\cdot a,\quad 
s^{(2)}_b(c):=c-I_a(b,c)\cdot b.
\end{eqnarray}
Then $s^{(1)}_a$ respects $I_s$ and is a reflection
(semisimple, eigenvalues $1,...,1,-1$).
And $s^{(2)}_b$ respects $I_a$ and is a pseudo-reflection
($s^{(2)}_b=\id$ or $s^{(2)}_b-\id$ nilpotent with one
single $2\times 2$ Jordan block).

(c) Now let $S=(s_{ij})\in T(n,\R)$
(so $s_{ij}=0$ for $i>j$, $s_{jj}=1$, 
and the eigenvalues of $S^{-1}S^t$ are in $S^1$). 
Then
\begin{eqnarray}\label{4.6}
(-1)^k\cdot M = s^{(k)}_{e_1}\circ ... \circ s^{(k)}_{e_n}
\qquad\textup{for }k\in\{1,2\}.
\end{eqnarray}
\end{theorem}

{\bf Proof:}
(a) $L(b,a)=L(Ma,b)$ is equivalent to 
$L(M\uuuu{e}^t,\uuuu{e}) = L(\uuuu{e}^t,\uuuu{e})^t$
which holds:
$$ L(M\uuuu{e}^t,\uuuu{e})
= L((\uuuu{e}\cdot S^{-1}S^t)^t,\uuuu{e}) 
=SS^{-t}\cdot S^t=S=L(\uuuu{e}^t,\uuuu{e})^t.$$
$M$ respects $L$ because of 
$$L(Ma,Mb)=L(Mb,a)=L(a,b).$$
$M$ respects $I_s$ and $I_a$ because of their relation
to $L$ in \eqref{4.2} and \eqref{4.3}.

\medskip
(b) $s^{(1)}_a$ respects $I_s$ because of 
\begin{eqnarray*}
&&I_s(s^{(1)}_a(b),s^{(1)}_a(c)) 
= I_s(b-I_s(a,b)a,c-I_s(a,c)a)\\
&=& I_s(b,c)-I_s(a,b)I_s(a,c) -I_s(a,c)I_s(b,a)
+I_s(a,b)I_s(a,c)I_s(a,a)\\
&=& I_s(b,c).
\end{eqnarray*}
$s^{(1)}_a$ is a reflection because its restriction to
$\{c\in H_\R\, |\, I_s(a,c)=0\}$ is $\id$ and because
of $s^{(1)}_a(a)=-a$. 

$s^{(2)}_b$ respects $I_a$ because of
\begin{eqnarray*}
&&I_a(s^{(2)}_b(c),s^{(2)}_b(d)) 
= I_a(c-I_a(b,c)b,d-I_a(b,d)b)\\
&=& I_a(c,d)-I_a(b,c)I_a(b,d) -I_a(b,d)I_a(c,b)
+I_a(b,c)I_a(b,d)I_a(b,b)\\
&=& I_a(c,d).
\end{eqnarray*}
$s^{(2)}_b$ is a pseudo-reflection because its restriction
to $\{c\in H_\R\, |\, I_a(b,c)=0\}$ is $\id$ and
this space has dimension $n-1$ or $n$ and contains $b$.

\medskip
(c) Denote $D_{kl}:= (\delta_{ik}\cdot \delta_{jl})_{i,j=1,...,n}
\in M(n\times n,\Z)$. Denote by 
$E_n:=(\delta_{ij})=\sum_{j=1}^n D_{jj}$
the $n\times n$ unit matrix.
Observe
\begin{eqnarray*}
D_{ij}D_{kl}=0\quad\textup{if }j\neq k,
\end{eqnarray*}
which implies 
\begin{eqnarray*}
(E_n+D_{ij})(E_n+D_{kl}) &=& E_n+D_{ij}+D_{kl}
\quad\textup{if }j\neq k,\\
(E_n+D_{ij})^{-1} &=& E_n-D_{ij}\quad\textup{if }i\neq j.
\end{eqnarray*}
These identities are applied often in the following
calculations. Empty places mean zeros.

{\tiny
\begin{eqnarray*}
S&=& 
\begin{pmatrix} 1 & & & \\ & \ddots & & \\ 
 & & 1 & s_{n-1,n} \\ & & & 1 \end{pmatrix}
\begin{pmatrix} 1 & & & & \\ 
 & \ddots & & & \\
 & & 1 & s_{n-2,n-1} & s_{n-2,n}\\ 
 & & & 1 & 0 \\ & & & & 1 \end{pmatrix}
...
\begin{pmatrix} 1 & s_{12} & ... & s_{1n} \\ & 1 & & \\ 
 & & \ddots &  \\ & & & 1 \end{pmatrix}   \\
&=& 
\begin{pmatrix} 1 & & & s_{1n} \\ & \ddots & & \vdots \\ 
 & & 1 & s_{n-1,n} \\ & & & 1 \end{pmatrix}
\begin{pmatrix} 1 & & & s_{1,n-1} & \\ 
 & \ddots & & \vdots & \\ 
 & & 1 & s_{n-2,n-1} \\ 
 & & & 1 & \\ 
 & & & & 1 \end{pmatrix}
...
\begin{pmatrix} 1 & s_{12} & &  \\ & 1 & & \\ 
 & & \ddots &  \\ & & & 1 \end{pmatrix} ,
\end{eqnarray*}
}

{\tiny
\begin{eqnarray*}
&&S^{-1}S^t \\
&=&
\begin{pmatrix} 1 & -s_{12} & ... & -s_{1n} \\
 & 1 & & \\  & & \ddots & \\ & & & 1 \end{pmatrix} ...
\begin{pmatrix} 1 & & & & \\ & \ddots & & & \\
 & & 1 & -s_{n-2,n-1} & -s_{n-2,n} \\
 & & & 1 & 0 \\ & & & & 1 \end{pmatrix}
\begin{pmatrix} 1 & & & \\ & \ddots & & \\ 
 & & 1 & -s_{n-1,n} \\ & & & 1 \end{pmatrix}\\
&& 
\begin{pmatrix} 1 & & & \\ s_{12} & 1 & & \\
 & & \ddots & \\ & & & 1 \end{pmatrix} ... 
\begin{pmatrix} 1 & & & & \\ & \ddots & & & \\ 
 & & 1 & & \\ s_{1,n-1} & ... & s_{n-2,n-1} & 1 & \\
 & & & & 1 \end{pmatrix}
\begin{pmatrix} 1 & & & \\ & \ddots & & \\ 
 & & 1 & \\ s_{1n} & ... & s_{n-1,n} & 1 \end{pmatrix}\\
&=&
\begin{pmatrix} 1 & -s_{12} & ... & -s_{1n} \\
 & 1 & & \\  & & \ddots  & \\ & & & 1 \end{pmatrix} 
\begin{pmatrix} 1 & & & & \\ 
s_{12} & 1 & -s_{23} & ... & -s_{2n} \\
 & & 1 & & \\ & & & \ddots &  \\ & & & & 1 \end{pmatrix} ... \\
&&
\begin{pmatrix} 1 & & & & \\ & \ddots & & & \\
 & & 1 & & \\
 s_{1,n-1} & ... & s_{n-2,n-1} & 1 & -s_{n-1,n} \\
 & & & & 1 \end{pmatrix}
\begin{pmatrix} 1 & & & \\ & \ddots & & \\ & & 1 & \\
 s_{1n} & ... & s_{n-1,n} & 1 \end{pmatrix}
\end{eqnarray*}
}
\begin{eqnarray*}
= (s^{(2)}_{e_1})^{mat}\cdot (s^{(2)}_{e_2})^{mat}\cdot ...
\cdot (s^{(2)}_{e_{n-1}})^{mat} \cdot (s^{(2)}_{e_n})^{mat}
\end{eqnarray*}
where the $n\times n$-matrix 
\begin{eqnarray*}
(s^{(2)}_{e_j})^{mat}:=
\begin{pmatrix} 1 & & & & & & \\ 
 & \ddots & & & & & \\ & & 1 & & & & \\
 s_{1j} & ... & s_{j-1,j} & 1 & -s_{j,j+1} & ... & -s_{jn} \\
 & & & & 1 & & \\ 
 & & & & & \ddots & \\ 
 & & & & & & 1 \end{pmatrix} 
\end{eqnarray*}
satisfies  
\begin{eqnarray*}
s^{(2)}_{e_j}\, \uuuu{e} = \uuuu{e} \cdot (s^{(2)}_{e_j})^{mat} .
\end{eqnarray*}
This shows $M=s^{(2)}_{e_1}\circ ...\circ s^{(2)}_{e_n}$.
Define the matrix $(s^{(1)}_{e_j})^{mat}$ by
\begin{eqnarray*}
s^{(1)}_{e_j}\, \uuuu{e} = \uuuu{e} \cdot (s^{(1)}_{e_j})^{mat} .
\end{eqnarray*}
Observe 
\begin{eqnarray*}
(s^{(1)}_{e_j})^{mat} &=&
\begin{pmatrix} 1 & & & & & & \\ 
 & \ddots & & & & & \\ & & 1 & & & & \\
 -s_{1j} & ... & -s_{j-1,j} & -1 & -s_{jj+1} & ... & -s_{jn} \\
 & & & & 1 & & \\ 
 & & & & & \ddots & \\ 
 & & & & & & 1 \end{pmatrix}  \\
&=& \Bigl(-\sum_{i=1}^{j-1}D_{ii}+\sum_{i=j}^nD_{ii} \Bigr)
\cdot (s^{(2)}_{e_j})^{mat} \cdot 
\Bigl(-\sum_{i=1}^{j}D_{ii}+\sum_{i=j+1}^nD_{ii} \Bigr)
\end{eqnarray*}
This shows 
\begin{eqnarray*}
-S^{-1}S^t = (s^{(1)}_{e_1})^{mat}\cdot ...\cdot
(s^{(1)}_{e_n})^{mat}\quad\textup{ and }\quad 
-M =s^{(1)}_{e_1}\circ ... \circ s^{(1)}_{e_n}.
\quad\Box
\end{eqnarray*}

\begin{corollary}\label{t4.2}
Consider the same situation as in theorem \ref{t4.1}.
Define the cyclic automorphism $C$ by 
\begin{eqnarray}\label{4.7}
C\, \uuuu{e} &=& \uuuu{e}\cdot 
\left(\begin{array}{ccc|c} & & & 1 \\ \hline 
& & & \\ & E_{n-1} & & \\ & & & \end{array}\right)
=\uuuu{e}\cdot C^{mat},\\
\textup{so }C \, e_j&=&e_{j+1}\textup{ for }1\leq j\leq n-1,
\ C \, e_n=e_1,\textup{ and }C^n=\id.\hspace*{1cm}\label{4.8}
\end{eqnarray}
Define the automorphisms $R_{(kj)}$ for
$k\in\{1,2\},j\in\{1,...,n\}$ of $H_\R$ by 
\begin{eqnarray}
R_{(kj)}&:=& C^{-(j-1)}\circ s^{(k)}_{e_j} \circ C^j.\label{4.9}
\end{eqnarray}
Then
\begin{eqnarray*}
R_{(kj)}\, \uuuu{e} = \uuuu{e}\cdot R_{(kj)}^{mat}
\end{eqnarray*}
with
\begin{eqnarray}\label{4.10}
R_{(1j)}^{mat} &=& 
\left(\begin{array}{cccccc|c} 
-s_{j,j+1} & ... & -s_{jn} & -s_{1j} & ... & -s_{j-1,j} & -1 \\
\hline
 & & & & & & \\
 & & & E_{n-1} & & & \\
 & & & & & & \end{array}\right),\hspace*{1cm}\\
R_{(2j)}^{mat} &=& 
\left(\begin{array}{cccccc|c} 
-s_{j,j+1} & ... & -s_{jn} & s_{1j} & ... & s_{j-1,j} & 1 \\
\hline
 & & & & & & \\
 & & & E_{n-1} & & & \\
 & & & & & & \end{array}\right), \label{4.11}
\end{eqnarray}
and
\begin{eqnarray}
(-1)^k\cdot S^{-1}S^t &=& R^{mat}_{(k1)}\cdot ...\cdot R^{mat}_{(kn)}
\nonumber\\
\textup{and}\quad (-1)^k\cdot M &=& R_{(k1)}\circ ...\circ R_{(kn)}.\label{4.12}
\end{eqnarray}
\end{corollary}

{\bf Proof:} 
$(-1)^kM=R_{(k1)}\circ ...\circ R_{(kn)}$ is an immediate 
consequence of 
$(-1)^kM =s^{(k)}_{e_1}\circ ...\circ s^{(k)}_{e_n}$
and the definition of $R_{(kj)}$ and $C^n=\id$.
The formulas for $R^{mat}_{(kj)}$ follow from the formulas
for $(s^{(k)}_{e_j})^{mat}$.\hfill$\Box$

\begin{remarks}\label{t4.3}
The matrices $R_{(kj)}^{mat}$ are {\it companion matrices}.
A companion matrix is here a matrix (empty places mean zeros)
\begin{eqnarray}\label{4.13}
R^{mat} &=& 
\left(\begin{array}{cccc|c} 
-p_{n-1} & -p_{n-2} & ... & -p_1 & -p_0 \\
\hline
 & & & &  \\
 & &  E_{n-1} & &  \\
 & & & &  \end{array}\right)
\end{eqnarray}
with $p_{n-1},...,p_0\in\C$. 
Its characteristic polynomial is $p(x)=x^n+p_{n-1}x^{n-1}
+...+p_1x+p_0$. 
For each eigenvalue $\kappa\in\C$, it has only one Jordan block.
A basis of a Jordan block of size $l+1$ with 
eigenvalue $\kappa$ is
\begin{eqnarray}\label{4.14}
v_j = \begin{pmatrix}
(n-1)_j\cdot\kappa^{n-1} \\
(n-2)_j\cdot\kappa^{n-2} \\
\vdots \\
(j)_j\cdot \kappa^j \\
0 \\ \vdots \\ 0 \end{pmatrix}
\quad\textup{for }j=0,1,...,l,
\end{eqnarray}
with
\begin{eqnarray}\label{4.15}
(a)_b:= a(a-1)\cdot ...\cdot (a-b+1)
\quad\textup{for }a\in\C,b\in\Z_{\geq 0},
\end{eqnarray}
(and $(a)_0=1$) and
\begin{eqnarray}\label{4.16}
\Bigl(\kappa^{-1}R^{mat}-E_n\Bigr)v_j
= j\cdot v_{j-1} \qquad (\textup{with }v_{-1}=0).
\end{eqnarray}
Here we used that $\kappa$ is a zero of 
$p^{(j)}(x)=(n)_jx^{n-j}+p_{n-1}(n-1)_jx^{n-1-j}+ ...
+p_j(j)_jx^0$ for $0\leq j\leq l$, and we used
\begin{eqnarray}\label{4.17}
(a)_b -(a-1)_b = b\cdot (a-1)_{b-1}\qquad 
\textup{for }b\in\Z_{\geq 1}.
\end{eqnarray}
\end{remarks}

\begin{definition}\label{t4.4}
Fix $n\in\Z_{\geq 1}$ and $k\in\{1,2\}$.

(a) The space of polynomials 
$T_{{\rm HOR}k}^{pol}(n,\R)\subset \R[x]_{\deg= n}$
was defined in definition \ref{t3.1}.
Define the map
\begin{eqnarray}\label{4.18}
S^{(k)}:T_{{\rm HOR}k}^{pol}(n,\R)&\to& GL(n,\R)\\
p(x)=x^n+p_{n-1}x^{n-1}+...+p_0&\mapsto &
\begin{pmatrix} 1 & p_{n-1} & ...& p_2 & p_1 \\
 & \ddots & \ddots & & p_2 \\ & & \ddots & \ddots & \vdots \\
 & & & \ddots & p_{n-1} \\ & & & & 1 \end{pmatrix}.\nonumber
\end{eqnarray}
Define its image as 
$T_{{\rm HOR}k}(n,\R):=S^{(k)}(T_{{\rm HOR}k}^{pol}(n,\R))$
(theorem \ref{t4.5} (a) will show that it is a subspace of
$T(n,\R)$). Define the map
\begin{eqnarray}\label{4.19}
R^{mat}_{(k)}: T_{{\rm HOR}k}(n,\R)&\to & GL(n,\R)\\
S=S^{(k)}(p)&\mapsto& \left(\begin{array}{cccc|c} 
-p_{n-1} & -p_{n-2} & ... & -p_1 & -p_0 \\
\hline
 & & & &  \\
 & &  E_{n-1} & &  \\
 & & & &  \end{array}\right) \nonumber
\end{eqnarray}
(recall $p_0=(-1)^{k-1}$). 
$R^{mat}_{(k)}(S)$ is a companion matrix, and
its characteristic polynomial is $p(x)$ 
by remark \ref{t4.3}.

\medskip
(b) For $S\in T_{{\rm HOR}k}(n,\R)$ take up the data in
theorem \ref{t4.1}. Define an automorphism 
$R_{(k)}(S):H_\R\to H_\R$ by 
$R_{(k)}(S)\, \uuuu{e}:= \uuuu{e}\cdot R^{mat}_{(k)}(S)$.

\medskip
(c) For $S\in T_{{\rm HOR}k}(n,\R)$ define 
$\Spp(S):=\Spp(p)$ and $\Sp(S):=\Sp(p)$ 
where $p\in T^{pol}_{{\rm HOR}k}(n,\R)$
is the characteristic polynomial of $R^{mat}_{(k)}(S)$ 
(or, equivalently, of $R_{(k)}(S)$), and where
$\Spp(p)$ and $\Sp(p)$ are defined in recipe \ref{t3.3} (c).
\end{definition}

Definition \ref{t4.4} (a) and the next formula \eqref{4.20}
are essentially due to Horocholyn \cite[ch. 2]{Ho17}
(he considered half of the cases).
He also studied the signature of $S+S^t$.
Theorem \ref{t4.5} and corollary \ref{t4.6} encompass his results.
In cases relevant for chain type singularities
(see section \ref{s7}), the matrices $S$ and $R^{mat}_{(k)}(S)$
are also given in \cite{OR77}.
But there \eqref{4.20} is not even mentioned, although
the authors are certainly aware of it.

\begin{theorem}\label{t4.5}
Choose $S\in T_{{\rm HOR}k}(n,\R)$ and take up the data in
theorem \ref{t4.1}.

(a) 
\begin{eqnarray}\label{4.20}
(-1)^k\cdot S^{-1}S^t = R^{mat}_{(k)}(S)^n\quad\textup{and}\quad 
(-1)^k\cdot M =R_{(k)}(S)^n.
\end{eqnarray}
The generalized eigenspaces of $R_{(k)}(S)$ are the spaces
$H_\kappa^{(R)}:=\ker((R_{(k)}(S)-\kappa\cdot\id)^n)\subset H_\C$
with $p(\kappa)=0$. 
The generalized eigenspaces of $M$ are the spaces
$H_\lambda =\bigoplus_{\kappa: (-1)^k\kappa^n=\lambda}H_\kappa^{(R)}$.
Especially, $T_{{\rm HOR}k}(n,\R)\subset T(n,\R)$. 
The monodromy $M$ and the automorphism $R_{(k)}(S)$
have a single Jordan block on $H_\kappa^{(R)}$  
(because of remark \ref{t4.3}).

\medskip
(b) $R_{(k)}(S)$ respects $L$. 
Therefore $H_\R$ decomposes $L$-orthogonally into
the Seifert form pairs $(H_1^{(R)}\cap H_\R,L)$,
$(H_{-1}^{(R)}\cap H_\R,L)$,
and $((H_\kappa^{(R)}\oplus H_{\oooo\kappa}^{(R)})\cap H_\R)$ 
for each $\kappa\in S^1$ with $\Imm\kappa>0$ and 
$H^{(R)}_\kappa\neq\{0\}$.

\medskip
(c) $\Spp(S)$ and the decomposition of $(H_\R,L)$ in (b) give a
polarized enhancement of $(H_\R,L)$ 
(definition \ref{t2.7}): $\Spp(S)$ consists of spp-ladders,
one for each eigenvalue $\kappa$ of $R_{(k)}(S)$.
The spp-ladder for $\kappa$ has length $l+1=\dim H^{(R)}_\kappa$, 
center $m=1$, and first spectral number
$\alpha$ with $e^{-2\pi i\alpha}=\kappa$
Furthermore
\begin{eqnarray}\label{4.21}
L(a,N^l\oooo{a})\in e^{\frac{1}{2}\pi i (2\alpha+l)}
\cdot \R_{>0}\quad
\textup{for }a\in H^{(R)}_\kappa - N(H^{(R)}_\kappa).
\end{eqnarray}
If $\kappa=\pm 1$, it is a single spp-ladder.
If $\kappa\neq\pm 1$, the partner spp-ladder 
is the one for $\oooo\kappa$.

\medskip
(d) The underlying spectrum $\Sp(S)$ is the one which
recipe \ref{t1.1} gives for $S$ if it is applied to
$T_{{\rm HOR}k}(n,\R)$ (see part (c) of theorem \ref{t1.3}).
\end{theorem}

{\bf Proof:}
(a) The coefficients $p_{n-1},...,p_1$ in the matrix 
\begin{eqnarray*}
S=\begin{pmatrix}1 & p_{n-1} & ... & p_1 \\ 
 & \ddots & \ddots & \vdots \\ 
 & & \ddots & p_{n-1} & \\ & & & 1 \end{pmatrix}
\in T_{{\rm HOR}k}(n,\R)
\end{eqnarray*} 
satisfy $p_{n-j}=(-1)^{k-1}p_j$. Therefore the matrices 
$R^{mat}_{(kj)}$ for $j\in\{1,...,n\}$ 
in corollary \ref{t4.2} are all equal
to one another and to $R^{mat}_{(k)}(S)$.
Thus $(-1)^k\cdot M=R_{(k)}(S)^n$ and \eqref{4.20}.
The other statements are immediate consequences of \eqref{4.20}.

(b) We have to prove 
$R^{mat}_{(k)}(S)^t\cdot S^t \cdot R^{mat}_{(k)}(S) = S^t$.
Equivalent is 
$S\cdot R^{mat}_{(k)}(S)=R^{mat}_{(k)}(S)^{-t}\cdot S$.
Recall $p_{n-j}=p_0\cdot p_j$ and observe
\begin{eqnarray*}
R^{mat}_{(k)}(S)^{-1}
&=& \left(\begin{array}{c|ccc}
 & & & \\ & & E_{n-1} & \\ & & &  \\ \hline 
 -p_0 & -p_1 & ... & -p_{n-1} \end{array}\right),\\
R^{mat}_{(k)}(S)^{-t} 
&=&\left(\begin{array}{ccc|c}
 & & & -p_0 \\ \hline & & &  -p_1 \\ & E_{n-1} & & \vdots \\
 & & & -p_{n-1} \end{array}\right).
\end{eqnarray*}
One calculates $S\cdot R^{mat}_{(k)}(S)$ and 
$R^{mat}_{(k)}(S)^{-t}\cdot S$ and finds in both cases
\begin{eqnarray*}
\left(\begin{array}{cccc|c}
0 & ... & ... & 0 & -p_0 \\ \hline 1 & p_{n-1} & ... & p_2 & 0 \\
 & \ddots & \ddots & \vdots & \vdots \\ 
 & & \ddots & p_{n-1} & \vdots \\
 & & & 1 & 0 \end{array}\right).
\end{eqnarray*}

(c) All statements in part (c) except that the enhancement is
polarized follow immediately from part (b) and from lemma \ref{t3.4} (b).

It rests to show that the enhancement is polarized, i.e. \eqref{4.21}.

$(-1)^k\cdot M=R_{(k)}(S)^n$ gives
$N=n\cdot (\textup{nilpotent part of }R_{(k)}(S))$.
On $H^{(R)}_\kappa$ 
\begin{eqnarray*}
N^l = n^l\cdot (\textup{nilpotent part of }R_{(k)}(S))^l
=n^l\cdot (\kappa^{-1} R_{(k)}(S)-\id)^l.
\end{eqnarray*}
The vector $v_l$ in remark \ref{t4.3} corresponds to an
element $a\in H^{(R)}_\kappa -N(H^{(R)}_\kappa)$.
We have to calculate the phase of
\begin{eqnarray*}
L(a,N^l\oooo{a})=
v_l^t\cdot S^t\cdot (\oooo{\kappa}^{-1}R^{mat}_{(k)}(S)-E_n)^l
\cdot \oooo{v_l}
= v_l^t\cdot S^t\cdot l!\cdot \oooo{v_0}
\end{eqnarray*}
and want to find $e^{\frac{1}{2}\pi i (2\alpha+l)}$.
We denote $p_n:=1$.

\begin{eqnarray*}
v_l^t\cdot S^t\cdot \oooo{v_0} 
= \begin{pmatrix}
(n-1)_l\kappa^{n-1} \\ (n-2)_l\kappa^{n-2} \\
\vdots \\ (l)_l\kappa^l \\ 0 \\ \vdots \\ 0 \end{pmatrix}^t
\begin{pmatrix} 1 & & & \\ p_{n-1} & \ddots & & & \\
 \vdots & \ddots & \ddots & \\ p_1 & ... & p_{n-1} & 1
\end{pmatrix}
\begin{pmatrix} \oooo{\kappa}^{n-1} \\ \oooo{\kappa}^{n-2} \\
\vdots \\ \oooo{\kappa}^0 \end{pmatrix} 
\end{eqnarray*}
\begin{eqnarray*}
&=& (n-1)_l\cdot \kappa^{n-1}\cdot \oooo{\kappa}^{n-1} \\
&+& (n-2)_l\cdot \kappa^{n-2}\cdot 
(p_{n-1}\cdot\oooo{\kappa}^{n-1}+p_n\cdot\oooo\kappa^{n-2}) \\
&+& ... \\
&+& (l)_l\cdot\kappa^l\cdot 
(p_{l+1}\cdot\oooo{\kappa}^{n-1}+ p_{l+2}\cdot\oooo{\kappa}^{n-2}
+ ... + p_n\cdot\oooo{\kappa}^l) \\
&=& ((n-1)_l+(n-2)_l+...+(l)_l)\cdot p_n\cdot\oooo{\kappa}^0 \\
&+& ((n-2)_l+...+(l)_l)\cdot p_{n-1}\cdot\oooo{\kappa}^1 
+ ... + (l)_l\cdot p_{l+1}\cdot \oooo{\kappa}^{n-l-1} \\
&=& \frac{1}{l+1}\Bigl[ 
(n)_{l+1}\cdot p_n\cdot\oooo{\kappa}^0 
+ (n-1)_{l+1}\cdot p_{n-1}\cdot\oooo{\kappa}^1 \\
&+& 
... + (l+1)_{l+1}\cdot p_{l+1}\cdot \oooo{\kappa}^{n-l-1}
\Bigr] \\
&=& \frac{1}{l+1}\cdot \oooo{\kappa}^{n-l-1}\cdot
p^{(l+1)}(\kappa).
\end{eqnarray*}
The last equality uses
\begin{eqnarray}\label{4.22}
(l+1)(n)_{l+1}=(n-1)_l+(n-2)_l+...+(l)_l,
\end{eqnarray}
which is an immediate consequence of \eqref{4.17}.

Now write $\uuuu{\beta}=(\beta_1,...,\beta_n)
:= \Pi^{-1}(p(x))$ and $\kappa_j:=e^{-2\pi i\beta_j}$.
Then $p(x) = \prod_{j=1}^n(x-\kappa_j)$ and 
$\kappa$ is a zero of it of order $l+1$.
Thus
\begin{eqnarray*}
p^{(l+1)}(\kappa) = (l+1)!\cdot \prod_{j:\, \kappa_j\neq \kappa}
(\kappa - \kappa_j).
\end{eqnarray*}
If $\kappa=\pm 1$ then a single spp-ladder
is associated to $H^{(R)}_\kappa$.
It satisfies $2\alpha+l=0$, so then 
\eqref{4.21} predicts $L(a,N^l\oooo{a})>0$,
so $v_l^t\cdot S^t\cdot \oooo{v_0}>0$.
Indeed, if $\kappa=1$ then the $\kappa_j\neq\kappa$
come in complex conjugate pairs  or are equal to $-1$,
so $p^{(l+1)}(1)>0$ and 
$v_l^t\cdot S^t \cdot \oooo{v_0}>0$.
If $\kappa=-1$ then the $\kappa_j\neq\kappa$ come 
in complex conjugate pairs or are equal to $1$.
Thus the multiplicity of $1$ is congruent to 
$n-l-1\mod 2$.
Therefore $p^{(l+1)}(-1)\in (-1)^{n-l-1}\cdot\R_{>0}$
and $v_l^t\cdot S^t\cdot \oooo{v_0}>0$. 

It rests to consider the case $\kappa\neq\pm 1$.
We can suppose $\Imm\kappa<0$.
Then an index $a$ exists with 
$\beta_{a-1}<\beta_a=...=\beta_{a+l}<\beta_{a+l+1}$
and $a+l\leq \frac{n}{2}$ and 
$\kappa = \kappa_a = ... = \kappa_{a+l}$. 

We have the four cases
$(k=1\, \&\, n\equiv 0(2))$, 
$(k=1\, \&\, n\equiv 1(2))$,
$(k=2\, \&\, n\equiv 0(2))$ and 
$(k=2\, \&\, n\equiv 1(2))$.
We treat only the case $(k=2\, \&\, n\equiv 0(2))$.
The other cases are analogous. Then $\kappa_1=1$, 
$\kappa_{\frac{n+2}{2}}=-1$ and 
\begin{eqnarray*}
&& \oooo{\kappa}^{n-l-1}\cdot\prod_{j:\, \kappa_j\neq\kappa}
(\kappa-\kappa_j) \\
&=& \oooo{\kappa}^{n-l-1}(\kappa-\kappa_1)
(\kappa-\kappa_{\frac{n+2}{2}})\cdot
\prod_{2\leq j\leq \frac{n}{2},\, \kappa_j\neq\kappa}
(\kappa^2-\kappa(\kappa_j+\oooo{\kappa_j})+1) \\
&=& \oooo{\kappa}^{n-l-1-(\frac{n-2}{2}-l-1)-1}
\cdot(\kappa-\oooo{\kappa})^{l+2}\cdot
\prod_{2\leq j\leq \frac{n}{2},\, \kappa_j\neq\kappa}
(\kappa+\oooo{\kappa}-(\kappa_j+\oooo{\kappa_j}))\\
&\in& \oooo{\kappa}^{n/2}\cdot (-i)^{l+2}\cdot (-1)^{a-2}
\cdot\R_{>0}.
\end{eqnarray*}
Here $\alpha=\alpha_{a+l}$ by the recipe \ref{t3.3}, and 
\begin{eqnarray*}
\oooo{\kappa}^{n/2}=(e^{2\pi i \beta_{a+l}})^{n/2}
=e^{\pi i n\cdot \beta_{a+l}}
=e^{\pi i(\alpha_{a+l}+n\cdot \gamma_{a+l})}
=e^{\pi i(\alpha+a+l-1)},\\
\oooo{\kappa}^{n/2}\cdot (-i)^{l+2}\cdot (-1)^{a-2}
= e^{\frac{1}{2}\pi i(2\alpha+l)}\cdot\R_{>0}.
\end{eqnarray*}

(d) This was essentially proved in the proof of theorem 
\ref{t1.3} (b). 
Define 
\begin{eqnarray*}
\uuuu\beta^{(k)}=(\beta_1^{(k)},...,\beta_n^{(k)})
:=(S^{(k)}\circ \Pi)^{-1}: 
T_{{\rm HOR}k}(n,\R)\to T^{scal}_{{\rm HOR}k}(n,\R).
\end{eqnarray*}
Then the functions $\beta_j^{(k)}:T_{{\rm HOR}k}(n,\R)\to [0,1]$ 
and the function $\alpha_k^{(k)}$ in the proof 
of theorem \ref{t1.3} (b) 
are related by the recipe \ref{t3.3} (a), i.e. by 
$\alpha_j^{(k)}=n\beta_j^{(k)}(S)-j+\frac{k}{2}$.
\hfill$\Box$

\bigskip

The following corollary of theorem \ref{t4.5}
gives an example, what is in the  
polarized enhancement in theorem \ref{t4.5} (c).
It was proved in a more elementary way in \cite{Ho17}
(for the cases considered there).

\begin{corollary}\label{t4.6}
Choose a matrix $S\in T_{{\rm HOR}k}(n,\R)$ and take
up the data in theorem \ref{t4.1}.
The symmetric form $I_s$ is nondegenerate on $H_{\neq -1}$.
Its signature on $H_\R\cap H_{\neq -1}$ is 
$(s_+,s_0,s_-)$ with  
\begin{eqnarray}\label{4.23}
s_+&=&|\{\alpha_j\, |\, \alpha_j\in(\frac{-1}{2},\frac{1}{2})
\mod 2\Z\}|,\\
s_-&=&\dim H_{\neq -1}-s_+,\quad s_0=0. \nonumber
\end{eqnarray}
\end{corollary}

{\bf Proof:} 
The polarized enhancement of $(H_\R,L)$ in theorem \ref{t4.5}
(c)  is (by remark \ref{t2.8}) a split Steenbrink 
polarized mixed Hodge structure on 
$H_\R\cong M(n\times 1,\R)$ of weight $m=1$. 
Such structures are studied in \cite{BH17}.
Theorem 4.6 in \cite{BH17} gives a square root of a Tate
twist, which allows to go from weight $m=1$ to an
arbitrary weight $\www m\in\Z$. 
In \cite[Corollary 3.13]{CKS86} (see also \cite[Theorem 7.5]{He03})
an equivalence between a polarized mixed Hodge structure
and a nilpotent orbit of polarized pure Hodge structures
is given. Especially, they have the same spectral numbers
and the same polarizing form. Therefore we can
work with a polarized pure Hodge structure of even weight 
$\www m$. In that case, the polarizing form on $H_{\neq -1}$ 
is $I_s$, and \eqref{4.23} is an immediate consequence of the polarization.
\hfill $\Box$

\begin{remark}\label{t4.7}
In corollary \ref{t4.6}, when is $I_s$ positive definite
on $H_\R$? Only if all spectral numbers are in
$(\frac{-1}{2},\frac{1}{2})\mod 2\Z$. But
by corollary \ref{t3.5} and remark \ref{t3.6}, 
the gaps between subsequent
spectral numbers (if they are ordered by size) are $\leq 1$.
This enforces that all spectral numbers are in 
$(\frac{-1}{2},\frac{1}{2})$. 
And this implies that the numbers $\beta_j$ in
$\uuuu\beta=(\Pi\circ S^{(k)})^{-1}(S)\in 
T^{scal}_{{\rm HOR}k}(n,\R)$ are {\it interlacing}
with the numbers $\gamma_1,...,\gamma_n$: 
Their pairwise distances are $|\beta_j-\gamma_j|<\frac{1}{2n}$.
Such an interlacing is also discussed in \cite{Ho17}.
\end{remark}

\begin{remark}\label{t4.8}
For $S\in T_{{\rm HOR}k}(n,\R)$ take up the data in theorem 
\ref{t4.1} and define $H_\Z:=M(n\times 1,\Z)$.
Then $L:H_\Z\times H_\Z\to\Z$ is unimodular, and
$R_{(k)}(S)$ and $M=(-1)^kR_{(k)}(S)^n$ are $L$-orthogonal
automorphisms of $H_\Z$. For $R_{(k)}(S)$ this follows from
theorem \ref{t4.5} (b).

Furthermore, let $\uuuu{e}^*$ be the $\Z$-basis of $H_\Z$ which
is left $L$-dual to the standard basis $\uuuu{e}$, i.e.
with $L((\uuuu{e}^*)^t,\uuuu{e})=E_n$. Then the matrix 
$R_{(k)}^{mat*}(S)$ of $R_{(k)}(S)$ with respect to $\uuuu{e}^*$,
so with $R_{(k)}(S)(\uuuu{e}^*) = \uuuu{e}^*\cdot 
R_{(k)}^{mat*}(S)$, is
\begin{eqnarray}\label{4.24}
R_{(k)}^{mat*}(S) = R_{(k)}^{mat}(S)^{-t}
= \left(\begin{array}{ccc|c}
 & & & -p_0 \\ \hline & & &  -p_1 \\ & E_{n-1} & & \vdots \\
 & & & -p_{n-1} \end{array}\right)
\end{eqnarray}
by the proof of theorem \ref{t4.5} (b).

This implies $R_{(k)}(S)(e_j^*)=e_{j+1}^*$ for $j\in\{1,...,n-1\}$.
So $R_{(k)}(S)$ is a cyclic automorphism of $H_\Z$.
This applies to the chain type singularities and is a remarkable
fact there (remark \ref{t7.4} (iv)).
\end{remark}

\section{The cases $n=2$ and $n=3$}\label{s5}
\setcounter{equation}{0}

\subsection{The case $n=2$}\label{s5.1}
Consider an upper triangular matrix 
$S=\begin{pmatrix}1&a\\0&1\end{pmatrix}$ with $a\in\R$
and consider the matrix 
$R^{mat}_{(1)}(S)=\begin{pmatrix}-a&-1\\1&0\end{pmatrix}$.
By the proof of theorem \ref{t4.5} (a) 
(or a direct calculation)
\begin{eqnarray}\label{5.1}
-S^{-1}S^t = R^{mat}_{(1)}(S)^2.
\end{eqnarray}
The characteristic polynomial of $R^{mat}_{(1)}(S)$ is
$p(x)=x^2+ax+1$. Thus $R^{mat}_{(1)}$ and
$S^{-1}S^t$ have eigenvalues in $S^1$ if and only if
$|a|\leq 2$. Therefore
\begin{eqnarray}\label{5.2}
T(2,\R)&=&T_{{\rm HOR}1}(2,\R) 
=\{\begin{pmatrix}1&a\\0&1\end{pmatrix}\, |\, 
a\in [-2,2]\}\cong [-2,2]\\
&\cong& T_{{\rm HOR}1}^{scal}(2,\R)
=\{(\beta_1,\beta_2)\, |\, \beta_1\in[0,\frac{1}{2}],
\beta_2=1-\beta_1\} \cong [0,\frac{1}{2}].\nonumber
\end{eqnarray}
The recipe \ref{t3.3} gives for $p(x)=x^2+ax+1$ with $|a|\leq 2$
\begin{eqnarray}\label{5.3}
\beta_1\in [0,\frac{1}{2}],\quad 
\beta_2=1-\beta_1\in [\frac{1}{2},1]\quad\textup{with }
2\cos(2\pi\beta_1)=-a,\\
\gamma_1=\frac{1}{4},\quad \gamma_2=\frac{3}{4},\label{5.4}\\
\alpha_1=2\beta_1-\frac{1}{2}\in [-\frac{1}{2},\frac{1}{2}],
\quad \alpha_2=2\beta_2-\frac{3}{2} =-\alpha_1.\label{5.5}
\end{eqnarray}

$\alpha_1$ is determined by $2\sin(\pi\alpha_1)=a$.
$R^{mat}_{(1)}(S)$ and $S^{-1}S^t$ are not semisimple 
precisely at the boundary of $T_{{\rm HOR}1}(2,\R)$.
There they have a $2\times 2$ Jordan block and the following
eigenvalues, and the spectral pairs are:
\begin{eqnarray}\label{5.6}
\begin{array}{r|l|l}
 & a=-2 & a=2 \\ \hline 
\beta_1 & 0 & \frac{1}{2} \\
\textup{eigenvalue of }R^{mat}_{(1)}(S) &
e^{-2\pi i \beta_1}=1 & e^{-2\pi i \beta_1}=-1 \\
\alpha_1 & -\frac{1}{2} & \frac{1}{2} \\
\textup{eigenvalue of }S^{-1}S^t &
 e^{-2\pi i \alpha_1}=-1 & e^{-2\pi i \alpha_1}=-1 \\
\Spp(S) & (-\frac{1}{2},2), (\frac{1}{2},0) &
(-\frac{1}{2},2), (\frac{1}{2},0) 
\end{array}
\end{eqnarray}
The following table lists the types of the Seifert form pairs
which one obtains by theorem \ref{t4.5} (c) for each
$a\in[-2,2]$.
\begin{eqnarray}\label{5.7}
\begin{array}{l|l}
a=0 & 2\cdot\Seif(1,1,1,1) \\ \hline 
a\in]-2,2[-\{0\} & 
\Seif(e^{-2\pi i\alpha_1},2,1,e^{\pi i\alpha_1}) \\ 
 & \cong \Seif(e^{2\pi i\alpha_1},2,1,e^{-\pi i\alpha_1})\\ \hline
a=\pm 2 & \Seif(-1,1,2,1)
\end{array}
\end{eqnarray}

The eigenvalue strata and the Seifert form strata 
(definition \ref{t1.5} (f)) in $T(2,\R)$ coincide.
One is $\{E_2\}$, the others are 
$\{\begin{pmatrix}1&a\\0&1\end{pmatrix},
\begin{pmatrix}1&-a\\0&1\end{pmatrix}\}$
for $a\in [-2,2]-\{0\}$. 

The set $T_{{HOR}2}(2,\R)$ has dimension $0$ by \eqref{1.2}.
It is $T_{{\rm HOR}2}(2,\R)=\{E_2\}$, and
\begin{eqnarray}\label{5.8}
R^{mat}_{(2)}(E_2) = \begin{pmatrix}0&1\\1&0\end{pmatrix},\quad
R^{mat}_{(2)}(E_2)^2 = E_2.
\end{eqnarray}
Recipe \ref{t3.3} gives in the case $k=2$ for $S=E_2$
\begin{eqnarray}\label{5.9}
\beta_1=0,\ \beta_2=\frac{1}{2},\ \gamma_1=0,\ 
\gamma_2=\frac{1}{2},\ \alpha_1=0,\ \alpha_2=0.
\end{eqnarray}

In the case $n=2$ conjecture \ref{t1.6} is satisfied 
(and conjecture \ref{t1.7} is empty). 
The only singularity up to suspension with $\mu=2$ is 
$A_2$. It is a chain type singularity. Theorem \ref{t7.6}
implies for $n=2$ conjecture \ref{t1.9} for function germs.

\subsection{The case $n=3$}\label{s5.2}
The following theorem \ref{t5.1} describes the set
$T(3,\R)$, its Seifert form strata and its eigenvalue strata
(definition \ref{t1.5} (f)). Define
\begin{eqnarray}\label{5.10}
&&f^\C:\C^3\to\C,\quad f(a_1,a_2,a_3):= 4+a_1a_2a_3-
(a_1^2+a_2^2+a_3^2),\hspace*{1.5cm}\\
&&f:= f^\C|_\R:\R^3\to \R,\nonumber
\end{eqnarray}
\begin{eqnarray}
S^{[3]}:\R^3&\to& M(3\times 3,\R),\nonumber\\
a=(a_1,a_2,a_3)&\mapsto& S^{[3]}(a)=
\begin{pmatrix}1&a_1&a_3\\ &1&a_2\\ &&1\end{pmatrix},\label{5.11}\\
M(3,\R)_{tri}
&:=& S^{[3]}(\R^3)\subset M(3\times 3,\R),\nonumber
\end{eqnarray}
\begin{eqnarray}\label{5.12}
\textup{ray}(S)&:=& S^{[3]}(\R\cdot a)\textup{ for }
S=S^{[3]}(a)\neq E_3 \ (\textup{i.e. for }a\neq 0).\hspace*{0.5cm}
\end{eqnarray}

\begin{theorem}\label{t5.1}
$T(3,\R)$ is the closed semialgebraic subset of $M(3,\R)_{tri}$
\begin{eqnarray}\label{5.13}
T(3,\R) &=& 
\{S^{[3]}(a)\in M(3,\R)_{tri}\, |\, 0\leq f(a_1,a_2,a_3)\leq 4\}.
\hspace*{0.5cm}
\end{eqnarray}
Consider the subsets
\begin{eqnarray}
T(3,\R)_{pos}&:=& \{S\in T(3,\R)\, |\, S+S^t 
\textup{ pos. def. or pos. semidefinite}\},\nonumber\\
T(3,\R)_{exc}&:=& \{S^{[3]}(2,2,2),S^{[3]}(-2,-2,2),\label{5.14} \\
&& \hspace*{0.5cm}S^{[3]}(-2,2,-2),S^{[3]}(2,-2,-2)\},\nonumber\\
T(3,\R)_{ind}&:=& T(3,\R)- T(3,\R)_{pos}.\nonumber
\end{eqnarray}
$T(3,\R)_{pos}$ is homeomorphic to a 3-ball and
$G_{\textup{sign},3}$-invariant ($G_{sign,n}$: definition \ref{t1.5} (e)).
\begin{eqnarray}\label{5.15}
\oooo{T(3,\R)_{ind}} = T(3,\R)_{ind}\cup T(3,\R)_{exc},\\
\oooo{T(3,\R)_{ind}}\cap T(3,\R)_{pos}=T(3,\R)_{exc}.
\nonumber
\end{eqnarray}
$\oooo{T(3,\R)_{ind}}$ is homeomorphic to four copies of
$[0,1]\times \R^2$. 
These components are permuted by the group $G_{\textup{sign},3}$.
Each component is in the open quadrant in 
$M(3,\R)_{tri}\cong\R^3$ which contains one of the points
in $T(3,\R)_{exc}$. 
The boundary $\partial T(3,\R)$ is smooth and transversal
to the rays $\textup{ray}(S)$ for 
$S\in M(3,\R)_{tri}-\{E_3\}$ except at the 4 points in
$T(3,\R)_{exc}$. At each of the 4 points in $T(3,\R)_{exc}$
it is isomorphic to a cone.

For each type of a Seifert form pair of rank 3, at most one 
Seifert form stratum exists. The following table lists
those which exist.

\begin{eqnarray*}
\begin{array}{l|l}
\textup{Type of a Seifert form pair} & 
\textup{description of Seifert form stratum} \\ \hline 
3\cdot\Seif(1,1,1,1) & \{E_3\} \\ \hline 
\Seif(1,1,1,1) &
\textup{diffeomorphic to a 2-sphere}\\ 
\hspace*{0.5cm} 
+\Seif(e^{-2\pi i \alpha_1},2,1,e^{\pi i\alpha_1}) &
\hspace*{0.5cm} \textup{in int}(T(3,\R)_{pos})  \\ \hline 
\Seif(1,1,1,1)+ & 
\partial T(3,\R)_{pos}-T(3,\R)_{exp} \\ 
\hspace*{0.5cm}  +\Seif(-1,1,2,1) & 
\hspace*{0.5cm} \approx \textup{2-sphere}-\textup{4 points}\\
\hline 
\Seif(1,1,1,1) & 
T(3,\R)_{exc} \\ 
\hspace*{0.5cm} + \Seif(-1,2,1) & \\ \hline 
\Seif(1,1,1,1) & 
\textup{the 4 components of }\partial T(3,\R)_{ind}
\textup{ whose}\\
\hspace*{0.5cm} 
+\Seif(-1,1,2,-1)  & 
\hspace*{0.5cm}\textup{closures contain points of }T(3,\R)_{exp} \\
 & \hspace*{0.5cm} \approx \textup{4 copies of }\R^2-\{0\} \\
 \hline 
\Seif(1,1,1,1) &
\textup{diffeomorphic to 4 copies of }\R^2, \\
\hspace*{0.5cm} 
+\Seif(e^{-2\pi i\alpha_1},2,1,-e^{\pi i\alpha_1}) & 
\hspace*{0.5cm}\textup{one in each component}\\
 & \hspace*{0.5cm} \textup{of int}(T(3,\R)_{ind}) \\ \hline 
\Seif(1,1,3,1) & 
\textup{the 4 components of }\partial\oooo{T(3,\R)_{ind}}\\
 & \hspace*{0.5cm}\textup{which do not intersect }T(3,\R)_{exc}\\ 
 & \hspace*{0.5cm}\approx\textup{4 copies of }\R^2
\end{array}
\end{eqnarray*}

The three Seifert form strata with eigenvalues $(1,-1,-1)$
form one eigenvalue stratum.
It is one component of $\partial T(3,\R)$. 
The other Seifert form strata are eigenvalue strata.
The following is a rough picture of a part of $T(3,\R)$.
The thick line indicates $T_{{\rm HOR}1}(3,\R)$,
which will be discussed below.
\end{theorem}

\noindent
\includegraphics[width=0.95\textwidth]{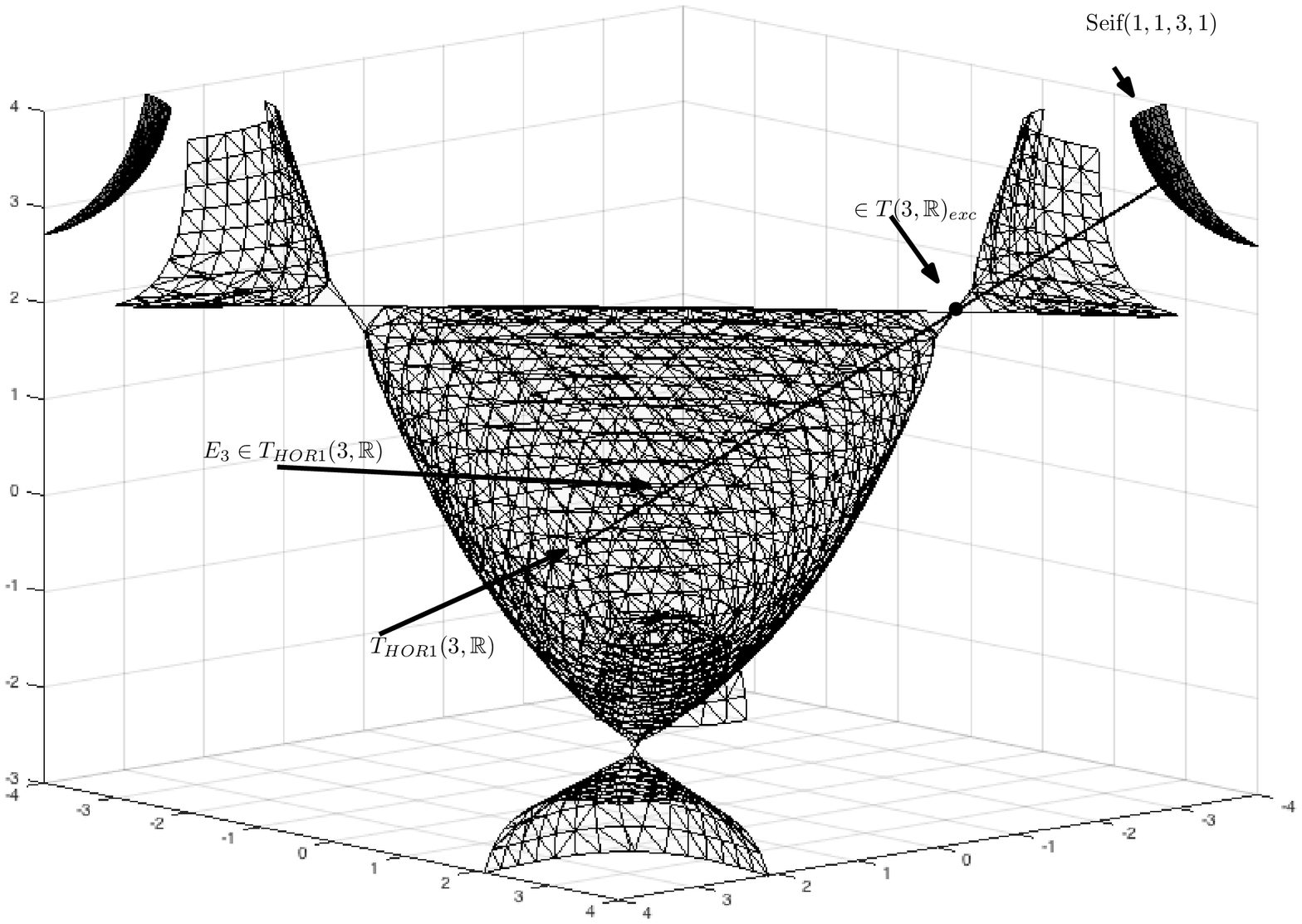}

{\bf Proof:}
The characteristic polynomial of $S^{-1}S^t$ for $S=S^{[3]}(a)$ is 
\begin{eqnarray}
p_{ch,S^{-1}S^t}(x) &=& \det(xE_3-S^{-1}S^t) = \det (xS-S^t) \nonumber\\
&=& x^3-(f(a)-1)x^2+(f(a)-1)x-1 \nonumber\\
&=& (x-1)(x^2-(f(a)-2)x+1) .\label{5.16}
\end{eqnarray}
This shows \eqref{5.13}.
The boundary $\partial T(3,\R)$ of $T(3,\R)$ is 
$\{S\in M(3,\R)_{tri}\, |\, f(a)=0\textup{ or }
f(a)=4\}-\{E_3\}$. 
For any $S=S^{[3]}(a)\in M(3,\R)_{tri}-\{E_3\}$,
consider the function
\begin{eqnarray*}
g^{\textup{ray},S}&:&\R_{\geq 0}\to \R,\\
g^{\textup{ray},S}(r)&:=& f(r\cdot a)
= 4+r^3\cdot a_1a_2a_3-r^2(a_1^2+a_2^2+a_3^2).
\end{eqnarray*}

\medskip
{\bf Claim 1:}
\begin{list}{}{}
\item[(i)]
If $a_1a_2a_3\leq 0$, then $g^{\textup{ray},S}$ is strictly 
decreasing with the limit $-\infty$,
so then $\textup{ray}(S)$ intersects $\partial T(3,\R)$
only in one point.
\item[(ii)]
If $a_1a_2a_3>0$ and $S\notin \textup{ray}(T(3,\R)_{exc})$,
then $g^{\textup{ray},S}$ is first strictly decreasing
to a minimum $<0$ and then strictly increasing with
limit $+\infty$.
Then $\textup{ray}(S)$ intersects $\partial T(3,\R)$
at three points.
\item[(iii)]
If $S\in T(3,\R)_{exc}$, then $g^{\textup{ray},S}$
is first strictly decreasing with minimum $=0$ at $S$
and then strictly increasing with limit $+\infty$.
Then $\textup{ray}(S)$ intersects $\partial T(3,\R)$
at $S$ and at one other point.
\end{list}

\medskip
{\bf Proof of claim 1:}
(i) is clear. (ii) and (iii): 
\begin{eqnarray*}
(g^{\textup{ray},S})'(r) &=& r\cdot (3r\cdot a_1a_2a_3
-2(a_1^2+a_2^2+a_3^2)),\\
r_0&:=& 2(a_1^2+a_2^2+a_3^2)/(3a_1a_2a_3),\quad
\textup{so that }(g^{\textup{ray},S})'(r_0)=0,\\
g^{\textup{ray},S}(r_0)&=& 
\frac{4}{27(a_1a_2a_3)^2}(27a_1^2a_2^2a_3^2 -
(a_1^2+a_2^2+a_3^2)^3)\\
&&\left\{ \begin{array}{ll}
=0 & \textup{ for }S\in T(3,\R)_{exc}, \\
\stackrel{(*)}{<} 0 &\textup{ for }a_1a_2a_3>0,
S\notin \textup{ray}(T(3,\R)_{exc}). \end{array} \right.
\end{eqnarray*}
$\stackrel{(*)}{<}$ is an easy exercise.
This finishes the proof of claim 1. \hfill $(\Box)$ 

\medskip
The eigenvalue map $\Psi_{Eig}:T(3,\R)\to\Eig(3)$
has the same fibers as the map $T(3,\R)\to\R,\ S(a)\to f(a)$.
Claim 1 shows that the fibers are smooth and transversal
to the rays $\textup{ray}(S)$, except at the point
$E_3$ and the four points in $T(3,\R)_{exc}$.
At each of these four points the fiber is locally diffeomorphic
to a cone, because $f$ has at
$(a_1,a_2,a_3)\in\{(2,2,2),(-2,-2,2),(-2,2,-2),(2,-2,-2)\}$
an $A_1$-singularity, and the signature of the Hessian
\begin{eqnarray*}
\textup{Hess}(f)(a) 
=\left(\frac{\paa^2 f}{\paa a_i\paa a_j}\right)(a)
=\begin{pmatrix}-2 & a_3 & a_2 \\ a_3 & -2 & a_1 \\
a_2 & a_1 & -2 \end{pmatrix}
\end{eqnarray*}
is $(1,0,2)$, because $\det \textup{Hess}(f)(a)=32>0$ 
and $-2<0$.

Claim 1 shows that $M(3,\R)_{tri}-\{S(a)\, |\, f(a)=0\}$
has six components:
the component $C_1$ which contains $E_3$,
the component $C_2$ which contains all of the four quadrants
with $a_1a_2a_3<0$ except their intersection with 
$\oooo{C_1}$,
and the four components $C_3,C_4,C_5,C_6$ which contain
each one of the partial rays in $S^{[3]}(\R_{>1}\cdot (S^{[3]})^{-1}(T(3,\R)_{exc}))$.

\eqref{5.16} implies 
\begin{eqnarray*}
\det (S+S^t) = 2\cdot f(a) 
\left\{\begin{array}{ll}
>0& \textup{ on }C_1,C_3,C_4,C_5\textup{ and }C_6,\\
<0& \textup{ on }C_2.\end{array}\right.
\end{eqnarray*}
$2\cdot E_3$ has signature $(3,0,0)$, any matrix $S+S^t$ with 
$S\in S^{[3]}(\R_{>1}\cdot (S^{[3]})^{-1}(T(3,\R)_{exc}))$
has signature $(1,0,2)$
because $\det\begin{pmatrix}2 & a_1\\ a_1&2\end{pmatrix}<0$
for such matrices. Therefore 
\begin{eqnarray*}
\textup{signature}(S+S^t)=\left\{\begin{array}{ll}
(3,0,0)& \textup{ on }C_1,\\
(2,0,1)& \textup{ on }C_2,\\
(1,0,2)& \textup{ on }C_3,C_4,C_5\textup{ and }C_6,
\end{array}\right.
\end{eqnarray*}
\begin{eqnarray*}
\textup{signature}(S+S^t)=\left\{\begin{array}{ll}
(2,1,0)& \textup{ on the part of }\paa T(3,\R) \\
 & \textup{ between }C_1\textup{ and }C_2,\\
(1,2,0)& \textup{ on }T(3,\R)_{exc},\\
(1,1,1)& \textup{ on the part of }\paa T(3,\R) \\
 & \textup{ between }C_2\textup{ and }C_3,C_4,C_5,C_6.
\end{array}\right.
\end{eqnarray*}
Thus $\dim\Rad(S+S^t)=1$ for $S$ in 
$\{S(a)\, |\, f(a)=0\}-T(3,\R)_{exc}$.
Therefore $S^{-1}S^t$ has for such an $S$ a
$2\times 2$ Jordan block with eigenvalues $-1$.
For $S\in T(3,\R)_{exc}$ it is semisimple with eigenvalues
$1,-1,-1$.

Finally, consider the set $\{S(a)\, |\, f(a)=4\}-\{E_3\}$.
It is the union of the four boundary components of $T(3,\R)$ 
which do not contain $T(3,\R)_{exc}$.
For $S\in \{S(a)\, |\, f(a)=4\}-\{E_3\}$, 
claim 1 gives $a_1a_2a_3>0$. This implies 
$\rk(S-S^t)=\rk\begin{pmatrix}0&a_1\\ -a_1&0\end{pmatrix}=2$
and $\dim\Rad(S-S^t)=1$ and that $S^{-1}S^t$ has
a single $3\times 3$ Jordan block with eigenvalue 1.

The proof up to now gives all statements in theorem \ref{t5.1}
except the table with Seifert form pairs and Seifert form strata.
The proof gives also the eigenvalue strata and the signature of $I_s$
at each point of $M(3,\R)_{tri}$. The table with Seifert form pairs
and Seifert form strata can now be deduced from the eigenvalues
of $S^{-1}S^t$, its Jordan block structure, claim 1,
the signature of $S+S^t$, and from lemma \ref{t2.4}.  \hfill$\Box$

\bigskip
Now we study the subvarieties 
$T_{{\rm HOR}k}(3,\R)\subset T(3,\R)$ for $k\in\{1,2\}$.
\begin{eqnarray*}
p(x) &=& x^3+p_2x^2+p_1x+p_0\\
&=&x^3+(-1)^{k-1}p_1x^2+p_1x+(-1)^{k-1}
\in T^{pol}_{{\rm HOR}k}(3,\R)\\
&=& \left\{\begin{array}{l}
x^3+p_1x^2+p_1x+1\\
\hspace*{1cm} =(x+1)(x^2+(p_1-1)x+1)\textup{ for }k=1,\\
x^3-p_1x^2+p_1x-1\\
\hspace*{1cm} =(x-1)(x^2-(p_1-1)x+1)\textup{ for }k=2.
\end{array}\right. 
\end{eqnarray*}
\begin{eqnarray*}
T_{{\rm HOR}1}(3,\R) &=& 
\begin{pmatrix}1&p_1&p_1\\ &1&p_1\\ &&1\end{pmatrix}\, |\, 
p_1\in[-1,3]\},\\
T_{{\rm HOR}2}(3,\R) &=& 
\begin{pmatrix}1&-p_1&p_1\\ &1&-p_1\\ &&1\end{pmatrix}\, |\, 
p_1\in[-1,3]\}.
\end{eqnarray*}
The element $g=(1,-1,1)\in G_{sign,3}$ exchanges 
$T_{{\rm HOR}1}(3,\R)$ and $T_{{\rm HOR}2}(3,\R)$.
In the picture after theorem \ref{t5.1}, the thick line indicates
$T_{{\rm HOR}1}(3,\R)$.

By lemma \ref{t3.4} (c) $\Spp(g(S))=\Spp(S)$ for 
$S\in \bigcup_{k=1,2}T_{{\rm HOR}k}(3,\R)$.
Therefore we restrict in the following to $T_{{\rm HOR}1}(3,\R)$.

\begin{theorem}\label{t5.2}
(a) $T_{{\rm HOR}1}(3,\R)$ intersects the Seifert form stratum
of type $\Seif(1,1,1,1)+\Seif(e^{-2\pi i \alpha_1},2,1,
e^{\pi i\alpha_1})$ twice, 
the Seifert form stratum of type
$\Seif(1,1,1,1)+\Seif(-1,1,2,-1)$ not at all and each other
Seifert form stratum once.

\medskip
(b) Recipe \ref{t3.3} gives for $T_{{\rm HOR}1}(3,\R)$
numbers $\beta_j,\gamma_j,\alpha_j$ for $j=1,2,3$ with
\begin{eqnarray}\label{5.17}
\beta_1\in[0,\frac{1}{2}],\ \beta_2=\frac{1}{2},\ 
\beta_3=1-\beta_1\in[\frac{1}{2},1],\\
\gamma_1 =\frac{1}{6},\ \gamma_2=\frac{1}{2},\ 
\gamma_3=\frac{5}{6},\nonumber\\
\alpha_1=3\beta_1-\frac{1}{2}\in [-\frac{1}{2},1],\ 
\alpha_2=3\beta_2-\frac{3}{2}=0,\nonumber\\
\alpha_3=3\beta_3-\frac{5}{2}=-\alpha_1\in[-1,\frac{1}{2}].
\nonumber
\end{eqnarray}
$\beta_1$ is determined by $\beta_1\in[0,\frac{1}{2}]$
and $\cos(2\pi \beta_1)=\frac{1-p_1}{2}$.
Thus $\beta_1$ and $\alpha_1$ are monotonically increasing with 
$p_1\in[-1,3]$.

\noindent
\includegraphics[width=1.0\textwidth]{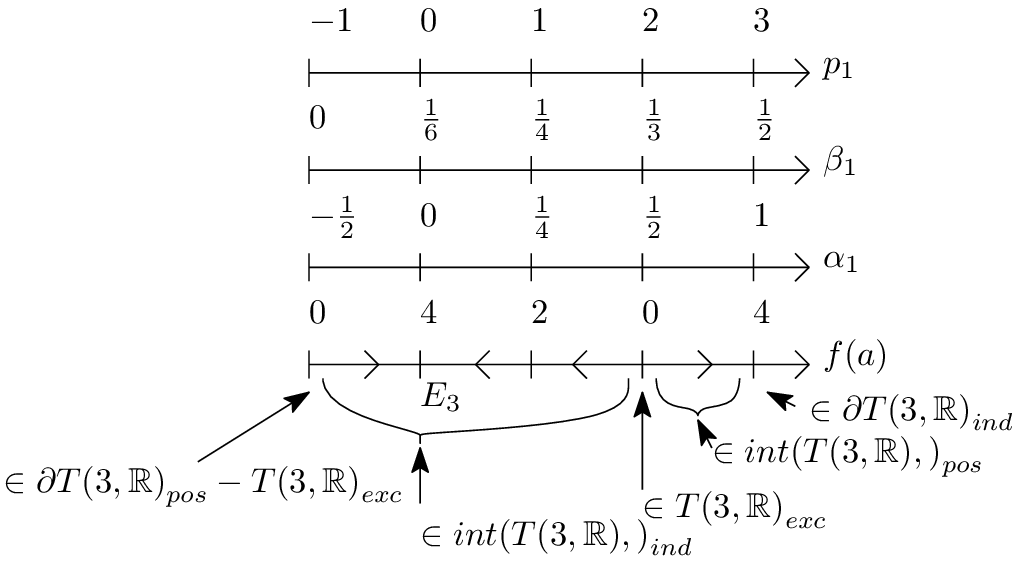}

(c) The conjectures \ref{t1.6} and \ref{t1.7} hold. The Seifert form
strata in $\oooo{T(3,\R)_{pos}}$ have spectral numbers
in $[-\frac{1}{2},\frac{1}{2}]$, 
the Seifert form strata in $\oooo{T(3,\R)_{ind}}$ 
have spectral numbers in $[-1,-\frac{1}{2}]
\cup \{0\}\cup[\frac{1}{2},1]$.
The two Seifert form strata with eigenvalues 
$1,-1,-1$ and a $2\times 2$ Jordan block for the eigenvalue
$-1$ have the same spectral pairs
$(0,1),(-\frac{1}{2},2),(\frac{1}{2},0)$.
The Seifert form stratum 
$\{S^{[3]}(a)\, |\, f(a)=4\}-\{E_3\}$
with a $3\times 3$ Jordan block has the spectral pairs
$(-1,3),(0,1),(1,-1)$.

\medskip
(d) Conjecture \ref{t1.9} for function germs holds in the case $n=3$.
\end{theorem}

{\bf Proof:}
(a) $T_{{\rm HOR}1}(3,\R)$ is the intersection of $T(3,\R)$
with the line through $E_3=S^{[3]}(0,0,0)$ and 
$S^{[3]}(2,2,2)\in T(3,\R)_{exc}$. This and theorem \ref{t5.1}
show part (a).

(b) $\beta_1$ in recipe \ref{t3.3} for 
$S\in T_{{\rm HOR}1}(3,\R)$ is determined by 
$\beta_1\in[0,\frac{1}{2}]$ and
$(x-e^{2\pi i\beta_1})(x-e^{-2\pi i\beta_1})
=x^2+(p_1-1)x+1$, which is
$\cos(2\pi \beta_1)=\frac{1-p_1}{2}$.
This shows all of (b).

(c) This follows from (a) and (b) and the following observation. 
At the boundary points of $T_{{\rm HOR}1}(3,\R)$,
the monodromy $S^{-1}S^t$ and $R^{mat}_{(1)}(S)$
have for each eigenvalue of
$R^{mat}_{(1)}(S)$ one Jordan block.

(d) The only singularity up to suspension with $\mu=3$ is
$A_3$. It is a chain type singularity.
Theorem \ref{t7.6} implies for $n=3$ conjecture \ref{t1.9}
for function germs. \hfill$\Box$

\begin{remarks}\label{t5.3}
(i) By theorem \ref{t5.2} (a), $T_{{\rm HOR}k}(3,\R)$
for $k\in\{1,2\}$ does not intersect the Seifert form
stratum of type $\Seif(1,1,1,1)+\Seif(-1,1,2,-1)$.
This is consistent with theorem \ref{t4.5} (c):
On this Seifert form stratum, the single spp-ladder
$(\frac{-1}{2},2),(\frac{1}{2},0)$ has first spectral
number $\alpha=\frac{-1}{2}$ and $l=1$, and
\begin{eqnarray*}
L(a,N^la)\in (-1)\cdot \R_{>0} 
=(-1)\cdot e^{\frac{1}{2}\pi i(2\alpha+l)}\cdot\R_{>0}.
\end{eqnarray*}
Theorem \ref{t4.5} (c) forbids the existence of a matrix
in $T_{{\rm HOR}k}(n,\R)$  and in this Seifert form stratum.

\medskip
(ii) The table \eqref{5.5} for $n=2$ and the table
in theorem \ref{t5.2} for $n=3$ show that precisely the
following  $S^1$-Seifert form pairs have no realization
as $(M(n\times 1,\R),L)$ with $L(a,b)=a^tS^tb$ and 
$S\in T(n,\R)$: 
All those for which $I_s$ is negative semidefinite
(cf. lemma \ref{t2.4} and remark \ref{t2.11} (vii)),
and all those which contain $\Seif(1,1,1,-1)$ or
$\Seif(1,1,3,-1)$. It is an interesting question
what holds for $n\geq 4$. 
\end{remarks}

\section{Isolated hypersurface singularities and 
$M$-tame functions}\label{s6}
\setcounter{equation}{0}

\noindent
The purpose of this section is merely to give references for
facts mentioned in the introduction, namely that holomorphic
functions germs $f:(\C^{m+1},0)\to(\C,0)$ with isolated
singularity at 0 (short: ihs)
and $M$-tame functions $f:X\to\C$ come equipped with
$\textup{Br}_\mu\ltimes G_{sign,\mu}$-orbits of 
{\it Stokes matrices} in $T(\mu,\Z)$ and with spectral pairs
$\Spp(f)$. First we recall the definition of an $M$-tame function.

\begin{definition}\cite{NZ90}\cite{NS99}
A function $f:X\to\C$ is $M$-tame if $X$ is an affine manifold
(of some dimension $m+1$) and if for some closed embedding
$X\hookrightarrow \C^N$ the following holds.
For any $\eta>0$ an $R(\eta)>0$ exists such that the fibers
$f^{-1}(\tau)$ with $|\tau|<\eta$ are transversal to all
spheres $S^{2N+1}_R=\{z\in \C^N\, |\, |z|=R\}$ with $R\geq R(\eta)$.
\end{definition}

Now we will treat the case 1: $f$ an ihs, and the case 2: 
$f$ M-tame, almost simultaneously. 
In both cases the definition domain shall have dimension $m+1$. 

In case 1, see e.g. \cite{Lo84}, \cite{AGV88} or \cite{Eb07} 
for the construction of a good representative 
$f:Y\to\Delta_\eta$ where 
$\Delta_\eta:=\{\tau\in\C\, |\, |\tau|<\eta\}$ is a sufficiently
small disk. In case 2, one can similarly construct a 
good representative $f:Y\to\Delta_\eta$ for $\eta>0$ sufficiently
large. The Milnor number $\mu$ is in case 1 the Milnor number
of the singularity at $0\in Y$ which is then the only singularity
of $f:Y\to\Delta_\eta$. In case 2, $\mu$ is the sum of the Milnor
numbers of all singularities of $f:Y\to\Delta_\eta$,
which are all singularities of $f:X\to\C$. 

In both cases, the relative homology groups (reduced if $m=0$)
$Ml(f,\zeta):=H_{m+1}(Y,f^{-1}(\zeta\eta),\Z)$ with $\zeta\in S^1$ 
are isomorphic  to $\Z^\mu$ \cite[(5.11)]{Lo84} \cite[ch. 2]{AGV88}, 
and some generators of them can be called (classes of)
{\it Lefschetz thimbles}. They form a flat $\Z$-lattice bundle on $S^1$. 
An intersection form for Lefschetz thimbles
is well defined on relative homology groups with different 
boundary parts. It is for any $\zeta\in S^1$ a $(-1)^{m+1}$
symmetric unimodular bilinear form
\begin{eqnarray}\label{6.1}
I_{Lef}:Ml(f,\zeta)\times Ml(f,-\zeta)\to\Z
\end{eqnarray}
Let $\gamma_\pi$ (respectively $\gamma_{-\pi}$) be the
isomorphism $Ml(f,-\zeta)\to Ml(f,\zeta)$ by flat shift in
mathematically positive (respectively negative) direction.
Then the classical Seifert form is given by
\begin{eqnarray}\label{6.2}
L: Ml(f,\zeta)\times Ml(f,\zeta)\to\Z,\quad 
L(a,b):=(-1)^{m+1}I_{Lef}(a,\gamma_{-\pi} b).
\end{eqnarray}
The classical monodromy $M$ 
and the intersection form $I$ on $Ml(f,\zeta)$ are given by
\begin{eqnarray}\label{6.3}
L(Ma,b)&=&(-1)^{m+1}L(b,a),\\
I(a,b)&=& -L(a,b)+(-1)^{m+1}L(b,a)=L((M-\id)a,b).\label{6.4}
\end{eqnarray}
We define a {\it normalized} Seifert form 
$L^{hnor}$ and a {\it normalized} monodromy $M^{hnor}$ by
\begin{eqnarray}\label{6.5}
L^{hnor}&:=&(-1)^{(m+1)(m+2)/2}\cdot L,\\
M^{hnor}&:=& (-1)^{m+1}M.\label{6.6}
\end{eqnarray}
Thus $M^{hnor}$ is the monodromy of $L$ and of $L^{hnor}$
in the sense of definition \ref{t2.2} (b).

Finally, we refer to \cite[ch. 2]{AGV88} or \cite[5.5]{Eb07}
for the definition of a distinguished basis 
$\uuuu\delta=(\delta_1,...,\delta_\mu)$ of $Ml(f,\zeta)$.
The set of distinguished bases forms one orbit of the group
$\textup{Br}_\mu\ltimes G_{sign,\mu}$.
Here $\textup{Br}_\mu$ is the braid group with $\mu$ strings,
and $G_{sign,\mu}$ was defined in definition \ref{t1.5} (e).
See \cite{AGV88} or \cite{Eb07} for the action of
$\textup{Br}_\mu$. The group $G_{sign,\mu}$ acts
componentwise by sign changes.
Each distinguished  basis $\uuuu\delta$ gives rise to 
one matrix 
\begin{eqnarray}\label{6.7}
S&:=& L^{hnor}(\uuuu\delta^t,\uuuu\delta)^t\in T(\mu,\Z).
\end{eqnarray}
We call these matrices {\it Stokes matrices} because
some of them encode certain Stokes structures
(which will not be discussed here).
These matrices form also one 
$\textup{Br}_\mu\ltimes G_{sign,\mu}$-orbit. 
In the case of the ihs, this orbit is finite only for the
simple and the simple elliptic singularities,
and the orbit of distinguished bases is finite only for
the simple singularities \cite{Eb16}.

\bigskip
Now we come to the spectral pairs.
In the case of an ihs $f$, spectral pairs $\Spp(f)$
were first defined by Steenbrink \cite{St77}
as invariants of his natural mixed Hodge
structure on the space dual to $Ml(f,1)$
(see also \cite{AGV88}).
It is in the notation of \cite{BH17} a signed 
Steenbrink polarized mixed Hodge structure of weight $m$.
For the polarization see \cite{He02} or \cite{BH17}.

In the case of an M-tame function $f$, the spectral pairs
are defined in the same way as invariants of Sabbah's
natural mixed Hodge structure \cite{Sa98} on the space
dual to $Ml(f,\zeta)$. A certain twist of Sabbah's Hodge
filtration is a part of a Steenbrink polarized
mixed Hodge structure of weight $m$ 
\cite[Corollary 11.4]{HS07} in the notation of \cite{BH17}.

In both cases, $f$ ihs or $f$ M-tame function,
$\Spp(f)$ is a union of single spp-ladders and
sppl-pairs with center $m$
(as the spectral pairs of any Steenbrink mixed Hodge
structure in the sense of \cite{BH17}).

\section{Chain type singularities and their spectra}\label{s7}
\setcounter{equation}{0}

\noindent
We used the initials of Horocholyn, Orlik and Randell
in the name HOR-matrices because a good part of these matrices
turns up in the papers \cite{Ho17} and \cite{OR77}.
See the beginning of section \ref{s4} for \cite{Ho17}.
Orlik and Randell studied the chain type singularities
(definition \ref{t7.1} below). They conjectured that each
of them has a distinguished basis whose Stokes matrix
$S$ is a certan HOR-matrix $S$ \cite[Conjecture (4.1)]{OR77}
(=conjecture \ref{t7.3}).
This conjecture and theorem \ref{t4.5} (a) would imply
that the matrix of the monodromy for this distinguished
basis is $(R^{mat}_{(k)})^\mu$ with $k\equiv m(2)$ 
(remark \ref{t7.4} (iii)).
The main result theorem (2.11) in \cite{OR77}
says that the matrix of the monodromy for some basis of the
Milnor lattice is this matrix.

We will recall the definition of a chain type singularity
and the HOR-matrix of Orlik and Randell.
Then we will show that the spectrum $\Sp(S)$
of this HOR-matrix from definition \ref{t4.4} (c)
(or theorem \ref{t1.3} (b), see theorem \ref{t4.5} (c)) 
is up to the shift
$\frac{m-1}{2}$ the correct spectrum of the singularity,
$\Sp(S)=\Sp(f)-\frac{m-1}{2}$.
This is positive evidence for conjecture \ref{t1.9}.
It is the main result of this section.
Of course, the evidence would be stronger if somebody
would prove conjecture (4.1) in \cite{OR77}.

\begin{definition}\label{t7.1}
(a) A {\it chain type singularity} is a function germ
on $(\C^{m+1},0)$ which is defined by a polynomial
\begin{eqnarray*}
f(x_0,...,x_m)=x_0^{a_0}+x_0x_1^{a_1}+...+x_{m-1}x_m^{a_m}
=x_0^{a_0}+\sum_{j=1}^m x_{j-1}x_j^{a_j}
\end{eqnarray*}
with $a_0\in \Z_{\geq 2}, a_1,...,a_m\in\Z_{\geq 1}$.

(b) Define the function
\begin{eqnarray}\label{7.1}
\rho:\bigcup_{k=0}^\infty \Z^k&\to&\Z,\\
\rho(a_0,a_1,...,a_{k-1})&:=& 
a_0...a_{k-1}-a_1...a_{k-1}+... +(-1)^{k-1}a_{k-1}
+(-1)^k\nonumber
\end{eqnarray}
(the case $k=0$ is $\rho(\emptyset)=1$).
\end{definition}

\begin{lemma}\label{t7.2}
Consider $f$ in definition \ref{t7.1} (a).
It has indeed an isolated singularity at 0.
It is a quasihomogeneous polynomial of weighted degree 1
with respect to weights $(w_0,...,w_m)$ which are determined
as follows. Define
\begin{eqnarray}\label{7.2}
r_{-1}:=1,\ r_k:=a_0...a_k=r_{k-1}a_k\quad\textup{for }
0\leq k\leq m,\\
\mu_{-1}:=1,\ \mu_k:=\rho(a_0,...,a_k)=r_k-\mu_{k-1}\quad
\textup{for }0\leq k\leq m,\label{7.3}\\
w_{-1}:=0,\ w_k:=\frac{\mu_{k-1}}{r_k}=\frac{1-w_{k-1}}{a_k}
\quad\textup{for }0\leq k\leq m.\label{7.4}
\end{eqnarray}
Its Milnor number is $\mu=\mu_m$.
\end{lemma}

{\bf Proof:}
The partial derivatives of $f$ are
\begin{eqnarray}\label{7.5}
\frac{\paa f}{\paa x_0} &=& 
a_0x_0^{a_0-1}+x_1^{a_1}, \\
\frac{\paa f}{\paa x_1} &=& 
a_1x_0x_1^{a_1-1}+x_2^{a_2},...,\nonumber\\
\frac{\paa f}{\paa x_{m-1}} &=& 
a_{m-1}x_{m-2}x_{m-1}^{a_{m-1}-1}+x_m^{a_m},\nonumber\\
\frac{\paa f}{\paa x_m} &=& 
a_mx_{m-1}x_m^{a_m-1}.\nonumber
\end{eqnarray}
Suppose that $x\in\C^{m+1}$ is a zero of all partial 
derivatives. Then
\begin{eqnarray*}
x_0\neq 0\Rightarrow x_1\neq 0\Rightarrow x_2\neq 0
\Rightarrow ...\Rightarrow x_m\neq 0\\
\Rightarrow \frac{\paa f}{\paa x_m}(x)\neq 0,
\textup{ a contradiction}.\\
x_0=0\Rightarrow x_1=0\Rightarrow x_2=0
\Rightarrow ...\Rightarrow x_m=0.
\end{eqnarray*}
Therefore the singularity $x=0$ of $f$ is the only singularity
in $\C^{m+1}$. The weights $(w_0,...,w_m)$
are uniquely determined by 
\begin{eqnarray*}
w_0&=&\frac{1}{a_0}=\frac{\mu_{-1}}{r_0},\\
w_k&=& \frac{1-w_{k-1}}{a_k}
=\frac{1-\frac{\mu_{k-2}}{r_{k-1}}}{a_k}
=\frac{r_{k-1}+\mu_{k-2}}{r_{k-1}a_k}
=\frac{\mu_{k-1}}{r_k}.
\end{eqnarray*}

In the following calculation of the Milnor number, 
$\stackrel{(*)}{=}$ is a well known formula for 
all quasihomogeneous singularities.
\begin{eqnarray*}
\mu \stackrel{(*)}{=} 
\prod_{k=0}^m \left(\frac{1}{w_k}-1\right)
=\prod_{k=0}^m\frac{r_k-\mu_{k-1}}{\mu_{k-1}}
=\prod_{k=0}^m \frac{\mu_k}{\mu_{k-1}} =\mu_m .
\hspace*{1cm}\Box
\end{eqnarray*}

\begin{conjecture}\label{t7.3}
\cite[Conjecture (4.1)]{OR77}
The chain type singularity 
$f=x_0^{a_0}+x_0x_1^{a_1}+...+x_{m-1}x_m^{a_m}$
has a distinguished basis whose Stokes matrix $S$
is the HOR-matrix $S$ (definition \ref{t4.4} (a))
with polynomial
\begin{eqnarray}\label{7.6}
p(x) = x^\mu +p_{\mu-1}x^{\mu-1}+...+p_0
=\prod_{k=-1}^m (x^{r_k}-1)^{(-1)^{m-k}}.
\end{eqnarray}
\end{conjecture}

\begin{remarks}\label{t7.4}
(i) In conjecture \ref{t7.3} $p(x)$ has only simple eigenvalues,
namely all zeros of $x^{r_m}-1$ minus certain gaps,
which are most zeros of $x^{r_{m-1}}-1$.

\medskip
(ii) In conjecture \ref{t7.3} $p_0=(-1)^{m+1}$
and $S\in T_{{\rm HOR},k}(\mu,\R)\cap T(\mu,\Z)$ with $k\equiv m(2)$.

\medskip
(iii) Theorem (2.11) in \cite{OR77} says that for a suitable
basis of $Ml(f)$, the monodromy matrix is
$R^{mat}_{(k)}(S)^{\mu}$ with $k\equiv m(2)$.
This is compatible with conjecture \ref{t7.3} and 
theorem \ref{t4.5} (a), which give this for a distinguished
basis with Stokes matrix $S$. Here recall that in the 
singularity case the monodromy in theorem \ref{t4.1}
is the normalized monodromy $M^{nor}$ in \eqref{6.6} 
and that the true  monodromy is $(-1)^{m+1} M^{nor}$.

\medskip
(iv) Conjecture \ref{t7.3} and definition \ref{t4.4} (b)
give the automorphism 
$R_{(k)}(S):Ml(f)\to Ml(f)$  (with $k\equiv m(2)$) 
with characteristic polynomial $p(x)$.
It respects $L$ by theorem \ref{t4.5} (b),
it satisfies $R_{(k)}(S)^\mu=M$ by theorem \ref{t4.5} (a),
and it is cyclic by remark \ref{t4.8}.
\end{remarks}

\begin{remarks}\label{t7.5}
Here we will argue that it is almost always
(and especially in the proof of theorem \ref{7.6})
sufficient to consider chain type singularities
$f=x_0^{a_0}+x_0x_1^{a_1}+...+x_{m-1}x_m^{a_m}$
with $a_0\in\Z_{\geq 3},a_1,...,a_m\in\Z_{\geq 2}$,
and the $A_1$-singularity $x_0^2$.

\medskip
(i) $f(x)$ is right equivalent to 
$c_0\cdot x_0^{a_0}+c_1\cdot x_0x_1^{a_1}+...+
c_m\cdot x_{m-1}x_m^{a_m}$ for arbitrary
$c_0,...,c_m\in\C^*$.

\medskip
(ii) Let $f(x_0,...,x_m)$ be a chain type singularity
with $a_0=2$. Consider the new coordinates
$\www{x}_0=x_0+\frac{1}{2}x_1^{a_1},
\www{x}_k=x_k$ for $1\leq k\leq m$. Then
\begin{eqnarray}
f(x_0,...,x_m)&=&
(x_0+\frac{1}{2}x_1^{a_1})^2 -\frac{1}{4}x_1^{2a_1}
+x_1x_2^{a_2}+...+x_{m-1}x_m^{a_m}\nonumber \\
&=& \www{x}_0^2 - \frac{1}{4}\www{x}_1^{2a_1}
+\www{x}_1\www{x}_2^{a_2}+...+\www{x}_{m-1}\www{x}_m^{a_m}.
\label{7.7}
\end{eqnarray}
This is (up to a rescaling in $\www{x}_1$)
a 1-fold suspension of the chain type singularity
$\www f(y_0,...,y_{m-1})
=y_0^{2a_1}+y_0y_1^{a_2}+...+y_{m-2}y_{m-1}^{a_m}$ with
\begin{eqnarray*}
\www a_0&=&2a_1,\quad \www a_k=a_{k+1}
\quad\textup{for }1\leq k\leq m-1,\\
\www r_k&=&r_{k+1}
\quad\textup{for }0\leq k\leq m-1,
\quad\textup{all }\www r_k\equiv 0(2),\\
p(x)&=& (x+1)^{(-1)^{m}}\cdot\prod_{k=1}^m 
(x^{r_k}-1)^{(-1)^{m-k}},\\
\www p(x)&=& (x-1)^{(-1)^{m}}\cdot\prod_{k=1}^m 
(x^{r_k}-1)^{(-1)^{m-k}}= (-1)^\mu\cdot p(-x),\\
\Sp(\www f)&=& \Sp(f)-\frac{1}{2}.
\end{eqnarray*}
Lemma \ref{t3.4} (c) implies 
$\Sp(\www p(x))=\Sp(p(x))$.

\medskip
(iii) Let $f(x_0,...,x_m)$ be a chain type singularity
with $a_0=3$. Suppose that it has an exponent
$a_j=1$ and that $a_1,...,a_{j-1}\geq 2$.
Consider the new coordinates 
$\www x_{j-1}=x_{j-1}+x_{j+1}^{a_{j+1}}$ and
$\www x_k=x_k$ for $k\neq j-1$. Then
\begin{eqnarray}
&&f(x_0,...,x_m)\nonumber \\
&=& x_0^{a_0}+x_0x_1^{a_1}+...+x_{j-2}x_{j-1}^{a_{j-1}}
+(x_{j-1}+x_{j+1}^{a_{j+1}})x_j \nonumber\\
&+& x_{j+1}x_{j+2}^{a_{j+2}}+...+x_{m-1}x_m^{a_m}\nonumber\\
&=& \www x_0^{a_0} + \www x_0\www x_1^{a_1} + ... + 
(-1)^{a_{j-1}}\www x_{j-2}\www x_{j+1}^{a_{j-1}a_{j+1}} + 
\www x_{j+1}\www x_{j+2}^{a_{j+2}} 
+ ... + \www x_{m-1} \www x_m^{a_m}\nonumber \\
&+& \www x_{j-1}\www x_j 
+ \www x_{j-2}\cdot\Bigl(
\sum_{k=0}^{a_{j-1}-1} (-1)^k 
\begin{pmatrix}a_{j-1}\\k\end{pmatrix}
\www x_{j-1}^{a_{j-1}-k} (\www x_{j+1}^{a_{j+1}})^k\Bigr).
\label{7.8}
\end{eqnarray}
The first line of \eqref{7.8} is a chain type singularity
$\www f(y_0,...,y_{m-2})$ with
\begin{eqnarray*}
\www a_k &=& a_k\textup{ for }0\leq k\leq j-2,
\www a_{j-1}=a_{j-1}a_{j+1},\\ 
&& \www a_k=a_{k+2}\textup{ for }j\leq k\leq m-2,\\
\www r_k &=& r_k\textup{ for }0\leq k\leq j-2,\ 
\www r_k=r_{k+2}\textup{ for }j-1\leq k\leq m-2,\\
\www p(x)&=& p(x).
\end{eqnarray*}
The first monomial $\www x_{j-1}\www x_j$ in the second line
of \eqref{7.8} gives a 2-fold suspension of $\www f$.
The second part $\www x_{j-2}\cdot (...)$ consists
of monomials of weighted degree $>1$
if one associates to $\www x_{j-1}$ and to $\www x_j$ the 
degree $\frac{1}{2}$, because for $k=a_{j-1}-1$
\begin{eqnarray*}
&& a_{j+1}a_{j-1}\www w_{j+1} = 1-\www w_{j-2}\textup{ and }\\
&&\www w_{j-2}+\www w_{j-1}+a_{j+1}(a_{j-1}-1)\www w_{j+1}\\
&=&\frac{1}{2}+1-a_{j+1}\www w_{j+1}
= \frac{1}{2}+1-\frac{1-\www w_{j-2}}{a_{j-1}}>1
\quad\textup{because }a_{j-1}\geq 2.
\end{eqnarray*}
Therefore $f$ is right equivalent to a 2-fold suspension of 
$\www f$. This implies
\begin{eqnarray*}
\Sp(\www f)&=& \Sp(f)-1.
\end{eqnarray*}

(iv) One transforms a chain type singularity with $a_0=2$
with (ii) to a 1-fold suspension of a chain type singularity with
one variable less. One repeats (ii) until one arrives either
at the $A_1$-singularity $x_0^2$ or at a chain type singularity
with $a_0\geq 3$. Then one repeats (iii) until one arrives
at a chain type singularity with $a_0\geq 3$,
$a_1,...,a_{\www m}\geq 2$. 
Then $\Sp(\www p(x))=\Sp(p(x))$.
\end{remarks}

\begin{theorem}\label{t7.6}
Consider a chain type singularity
$f(x)=x_0^{a_0}+x_0x_1^{a_1}+...+x_{m-1}x_m^{a_m}$.
The spectrum of the HOR-matrix $S$ in conjecture \ref{t7.3}
(see definition \ref{t4.4} (c) for $\Sp(S)$) 
satisfies
\begin{eqnarray}\label{7.9}
\Sp(S)=\Sp(f)-\frac{m-1}{2}.
\end{eqnarray}
\end{theorem}

{\bf Proof:}
For the $A_1$-singularity $x_0^2$
$S=(1)$ and $\Sp(S)=(0)$ and $\Sp(f)=(-\frac{1}{2})$
and $m=0$, so \eqref{7.9} holds.
Because of this and the remarks \ref{t7.5} (ii)--(iv),
it is sufficient to prove theorem \ref{t7.6} for the cases
$a_0\in \Z_{\geq 3},a_1,...,a_m\in\Z_{\geq 2}$.
The spectrum $\Sp(f)=(\alpha_1(f),...,\alpha_\mu(f))$
(with an arbitrary numbering) of a quasihomogeneous
singularity with weights $w_0,...,w_m\in\Q\cap (0,1)$
such that $\deg_wf=1$ can be given in several ways:

\begin{list}{}{}
\item[(A)]
By the generating function
\begin{eqnarray}\label{7.10}
\sum_{j=1}^\mu t^{\alpha_j(f)+1} = \prod_{k=0}^m
\frac{t-t^{w_k}}{t^{w_k}-1}.
\end{eqnarray}
\item[(B)]
If $m_1,...,m_\mu\in\C[x]$ are weighted homogeneous 
polynomials which represent a basis of the Jacobi algebra then
\begin{eqnarray}\label{7.11}
\alpha_j(f)= -1+\sum_{k=0}^m w_k + \deg_w m_j
\quad\textup{for }j=1,...,\mu.
\end{eqnarray}
\end{list}

Here (B) is more convenient than (A).
Claim 1 is the first of four steps of the main part of the proof.

\medskip
{\bf Step 1 = Claim 1:}
The following monomials represent a basis of the Jacobi algebra:
\begin{eqnarray}
x_0^{b_0}x_1^{b_1}\cdot ...\cdot x_m^{b_m}
&\textup{with}& \left\{\begin{array}{l}
0\leq b_j\leq a_j-1\textup{ for }\\
\hspace*{1cm}j\in\{0,1,...,m-1\}\nonumber\\
\textup{and }0\leq b_m\leq a_m-2,\end{array}\right. \\
x_0^{b_0}x_1^{b_1}\cdot ...\cdot x_{m-2}^{b_{m-2}}x_m^{a_m-1}
&\textup{with}& \left\{\begin{array}{l}
0\leq b_j\leq a_j-1\textup{ for }\\
\hspace*{1cm}j\in\{0,1,...,m-3\}\\
\textup{and }0\leq b_{m-2}\leq a_{m-2}-2,\end{array}\right. 
\nonumber\\
\vdots && \nonumber\\
\textup{for }m\equiv 0(2):&&\nonumber\\
x_0^{b_0}x_2^{a_2-1}x_4^{a_4-1}\cdot ...\cdot x_m^{a_m-1} 
&\textup{with}&0\leq b_0\leq a_0-2, \nonumber\\
\textup{for }m\equiv 1(2):\nonumber\\
x_0^{b_0}x_1^{b_1}x_3^{a_3-1}\cdot ...\cdot x_m^{a_m-1} 
&\textup{with}& \left\{\begin{array}{l}
0\leq b_0\leq a_0-1\\ \textup{and }0\leq b_1\leq a_1-2 ,
\end{array}\right. \nonumber\\
x_1^{a_1-1}x_3^{a_3-1}\cdot ...\cdot x_m^{a_m-1}.&&\label{7.12}
\end{eqnarray}

\medskip
{\bf Proof of claim 1:} Their number is
\begin{eqnarray*}
\mu &=& a_0...a_{m-1}(a_m-1) + a_0...a_{m-3}(a_{m-2}-1) \\
&& + ... + \left\{\begin{array}{ll}
a_0-1 & \textup{ for }m\equiv 0(2)\\
a_0(a_1-1)+1 & \textup{ for }m\equiv 1(2)\end{array}\right. 
\end{eqnarray*}

Therefore for claim 1 it is sufficient to prove that any monomial
in $\C\{x\}$ is a linear combination of the monomials above
and of an element of the Jacobi ideal
$J_f =\left(\frac{\paa f}{\paa x_0},...,\frac{\paa f}{\paa x_m}
\right)$. The generators $\frac{\paa f}{\paa x_j}$ of $J_f$
are given in \eqref{7.5}. Obviously also
\begin{eqnarray*}
x_{m-1}x_m^{a_m},\ x_{m-2}x_{m-1}^{a_{m-1}},...,x_1x_2^{a_2},\ 
x_0x_1^{a_1},\ x_0^{a_0}
\end{eqnarray*}
are in $J_f$. Start with any monomial in $\C\{x\}$. Using
$\frac{\paa f}{\paa x_{m-1}}$, $\frac{\paa f}{\paa x_{m-2}}$,
...,$\frac{\paa f}{\paa x_0}$ and $x_0^{a_0}$
(in this order), one can reduce it modulo $J_f$ to 0 or
to a monomial $x_0^{b_0}\cdot ...\cdot x_m^{b_m}$
with $0\leq b_j\leq a_j-1$ for all $j$.

If $b_m\leq a_m-2$ stop here.
Suppose $b_m=a_m-1$. If $b_{m-1}\geq 1$ the monomial is in $J_f$.
Suppose $b_{m-1}=0$. If $b_{m-2}\leq a_{m-2}-2$ stop here.
Suppose $b_{m-2}=a_{m-2}-1$.
If $b_{m-3}\geq 1$, the monomial is modulo
$\C\cdot\frac{f}{\paa x_{m-2}}$ congruent to a monomial
$x_0^{\www b_0}\cdot ...\cdot x_m^{\www b_m}$ with 
$b_{m-1}\geq a_{m-1},\ b_m=a_m-1,$ so it is in $J_f$.
Suppose $b_{m-3}=0$. The claim is proved by repeating
these arguments. \hfill ($\Box$)

\medskip
{\bf Step 2:} 
The second step is the definition of a directed graph $G$
whose vertices are labelled by the monomials in \eqref{7.12}.
Before defining the directed edges, consider the following
$m+1$ Laurent monomials in $\C[x_0^{\pm 1},...,x_m^{\pm 1}]$:
\begin{eqnarray*}
\begin{array}{lcr}
\uuuu{x}^{\uuuu{g}(m)}&:=&x_m^{-1}, \\
\uuuu{x}^{\uuuu{g}(m-1)}&:=&x_{m-1}x_m^{a_m-2},\\ 
\uuuu{x}^{\uuuu{g}(m-2)}&:=&x_{m-2}^{-1}x_{m-1}^{-(a_{m-1}-1)}
x_m^{a_m-1}, \\
\uuuu{x}^{\uuuu{g}(m-3)}&:=&x_{m-3}^1x_{m-2}^{a_{m-2}-1}
x_{m-1}^{-(a_{m-1}-1)}x_m^{a_m-2},\\ 
\uuuu{x}^{\uuuu{g}(m-4)}&:=&x_{m-4}^{-1}x_{m-3}^{-(a_{m-3}-1)}
x_{m-2}^{a_{m-2}-1}x_{m-1}^{-(a_{m-1}-1)}x_m^{a_m-1}, \\
\uuuu{x}^{\uuuu{g}(m-5)}&:=&x_{m-5}^{1}x_{m-4}^{a_{m-4}-1}
x_{m-3}^{-(a_{m-3}-1)}x_{m-2}^{a_{m-2}-1}
x_{m-1}^{-(a_{m-1}-1)}x_m^{a_m-2},\\
&\vdots&\\
&& \textup{for }m\equiv 0(2): \hspace*{6cm}\\
\uuuu{x}^{\uuuu{g}(0)}&:=&x_0^{-1}x_1^{-(a_1-1)}...
x_{m-4}^{a_{m-4}-1}x_{m-3}^{-(a_{m-3}-1)}
x_{m-2}^{a_{m-2}-1}x_{m-1}^{-(a_{m-1}-1)}x_m^{a_m-1},\\
&&\textup{for }m\equiv 1(2):  \hspace*{6cm}\\
\uuuu{x}^{\uuuu{g}(0)}&:=&x_0^{1}x_1^{a_1-1}...
x_{m-4}^{a_{m-4}-1}x_{m-3}^{-(a_{m-3}-1)}
x_{m-2}^{a_{m-2}-1}x_{m-1}^{-(a_{m-1}-1)}x_m^{a_m-2}.\end{array}
\end{eqnarray*}
Now an edge labelled by $\uuuu g(j)$ goes from 
$\uuuu x^{\uuuu b}=x_0^{b_0}\cdot ...\cdot x_m^{b_m}$ 
to $\uuuu x^{\uuuu c}=x_0^{c_0}\cdot ...\cdot x_m^{c_m}$ if 
$\uuuu{x}^{\uuuu b}\cdot\uuuu{x}^{\uuuu{g}(j)} 
=\uuuu{x}^{\uuuu c}$. This defines a directed
graph $G$ with vertices labelled by the monomials
in \eqref{7.12} and edges labelled by $\uuuu g(0),...,
\uuuu g(m)$.

\medskip
{\bf Claim 2:}
(a) The graph is a chain. If $m\equiv 0(2)$ it starts at
$x_0^{a_0-1}x_2^{a_2-1}\cdot ...\cdot x_m^{a_m-2}$
and ends at
$x_1^{a_1-1}x_3^{a_3-1}\cdot ...\cdot x_{m-1}^{a_{m-1}-1}$. 
If $m\equiv 1(2)$ it starts at
$x_1^{a_1-1}x_3^{a_3-1}\cdot ...\cdot x_m^{a_m-1}$
and ends at 
$x_0^{a_0-1}x_2^{a_2-1}\cdot ...\cdot x_{m-1}^{a_{m-1}-1}$.

(b) The weight of the Laurent monomial
$\uuuu{x}^{\uuuu g(j)}$ is 
\begin{eqnarray}\label{7.14}
\deg_w \uuuu{x}^{\uuuu g(j)} = \left\{\begin{array}{ll}
-w_m & \textup{ if }j\equiv m(2), \\
1-2w_m & \textup{ if }j\equiv m+1(2).\end{array}\right.
\end{eqnarray}

\medskip
{\bf Proof of claim 2:}
(a) Careful inspection of the set of monomials in \eqref{7.12}.

(b) In both cases use $w_{k-1}+a_kw_k=1$. \hspace*{1cm}$\Box$

\medskip
{\bf Step 3:} 
The third step consists in making precise the recipe \ref{t3.3} in the 
case of the HOR-matrix $S$ respectively its polynomial
$p(x)$ in conjecture \ref{t7.3}.
Because $p_0=(-1)^{m+1}$, the case $m\equiv 0(2)$ is the
case $k=1$ in recipe \ref{t3.3}, and the case $m\equiv 1(2)$
is the case $k=2$ in recipe \ref{t3.3}. Then
\begin{eqnarray*}
\alpha_j &=& \mu(\beta_j-\gamma_j) \quad\textup{for }j=1,...,\mu\\
\mu\cdot\gamma_j &=& \left\{\begin{array}{ll}
j-\frac{1}{2} & \textup{ in the case }m\equiv 0(2)\\
j-1 & \textup{ in the case }m\equiv 1(2) \end{array}\right. \\
\mu\cdot \uuuu\beta &=& 
\frac{\mu}{r_m}(\delta_1,...,\delta_\mu)\quad
\textup{with }\delta_{j+1}=\delta_j+1\textup{ or }
\delta_{j+1}=\delta_j+2,\\
&& \{\delta_1,...,\delta_\mu\} \subset \{0,1,2,...,r_m-1\},
\textup{ namely}
\end{eqnarray*}
\begin{eqnarray*}
\prod_{j=1}^\mu (x-e^{-2\pi i \delta_j/r_m})
=p(x) = \prod_{l=0}^m (x^{r_l}-1)^{(-1)^{m-l}}.
\end{eqnarray*}
$\alpha_1,...,\alpha_\mu$ denote now the spectral numbers
in $\Sp(S)$ with the order from recipe \ref{t3.3}.
We find:
\begin{eqnarray}\label{7.15}
\begin{array}{lll}
\textup{If }\delta_{j+1}=\delta_j+1 & \textup{then} & 
\alpha_{j+1}-\alpha_j = \frac{\mu}{r_m}-1 
= \frac{\mu-r_m}{r_m} = \frac{ -\mu_{m-1}}{r_m} =-w_m \\
\textup{If }\delta_{j+1}=\delta_j+2 & \textup{then} & 
\alpha_{j+1}-\alpha_j = 2\frac{\mu}{r_m}-1 
= 1-2w_m .\end{array}
\end{eqnarray}

We have to show that $\alpha_1,...,\alpha_\mu$ coincide
up to the shift by $\frac{m-1}{2}$ with the spectral
numbers of $f$ which are given by \eqref{7.11} 
and \eqref{7.12}.

\medskip
{\bf Step 4 = Claim 3:}
Denote the monomials in \eqref{7.12} by $m_1,...,m_\mu$
with the numbering as the chain $G$ prescribes it.
Denote $\alpha_1(f),...,\alpha_\mu(f)$ according
to \eqref{7.11}. Then
\begin{eqnarray*}
\alpha_j = \alpha_j(f) -\frac{m-1}{2},\quad
\textup{so }\Sp(S)=\Sp(f)-\frac{m-1}{2}.
\end{eqnarray*}

\medskip
{\bf Proof of claim 3:}
If the vertices $m_j$ and $m_{j+1}$ in the chain $G$ are 
connected by an edge of type $\uuuu g(l)$ then
\begin{eqnarray*}
\alpha_{j+1}(f)-\alpha_j(f) 
&=& \deg_w m_{j+1}-\deg_w m_j 
=\deg_w \uuuu{x}^{\uuuu g(l)} \\
&=& \left\{\begin{array}{ll} 
-w_m & \textup{ if }l\equiv m(2) \\
1-2w_m & \textup{ if }l\equiv m+1(2) \end{array}\right. 
\end{eqnarray*}
Therefore it rests to see two points:
\begin{eqnarray*}
\alpha_1(f) &=& \frac{m-1}{2} +\alpha_1, \\
\delta_{j+1}=\delta_j+2 &\iff& 
\textup{the edge from }m_j \textup{ to }m_{j+1}
\textup{ is of type }\\
&&\gamma_l\textup{ with }l\equiv m+1(2).
\end{eqnarray*}
We carry out the first point in both cases $m\equiv 0(2)$
and $m\equiv 1(2)$ and leave the second point to the reader.

\medskip
The case $m\equiv 0(2):$  Then
\begin{eqnarray*}
\alpha_1 &=& \mu(\beta_1-\gamma_1)
=\frac{\mu}{r_m}-\frac{1}{2}=\frac{1}{2}-w_m,\\
\alpha_1(f) &=& -1+\sum_{k=0}^mw_k + 
\deg_w x_0^{a_0-1}x_2^{a_2-1}...x_m^{a_m-2} \\
&=& -1+\deg_w x_0^{a_0}x_1x_2^{a_2}x_3...x_{m-1}x_m^{a_m-1}\\
&=& \frac{m}{2}-w_m =\frac{m-1}{2}+\alpha_1.
\end{eqnarray*}

\medskip
The case $m\equiv 1(2):$  Then
\begin{eqnarray*}
\alpha_1 &=& \mu(\beta_1-\gamma_1)
=\frac{0}{r_m}-0=0,\\
\alpha_1(f) &=& -1+\sum_{k=0}^mw_k +  
\deg_w x_1^{a_1-1}x_3{a_3-1}...x_m^{a_m-1} \\
&=& -1 + 
\deg_w x_0x_1^{a_1}x_2x_3^{a_3}...
x_{m-2}^{a_{m-2}}x_{m-1}x_m^{a_m} \\
&=& \frac{m-1}{2}=\frac{m-1}{2}+\alpha_1.\hspace*{1cm}(\Box)
\end{eqnarray*}
This finishes the proof of theorem \ref{t7.6} \hfill$\Box$

\section{Some remarks and speculations}\label{s8}
\setcounter{equation}{0}

\noindent
In the following three subsections, we offer a critical
discussion of some arguments in \cite{CV93} with a counterexample,
we make a few comments on flat vector bundles and a few
comments on Thom-Sebastiani formulas.

\subsection{Arguments in \cite{CV93} for conjecture \ref{t1.9}}
\label{s8.1}

The arguments concern the case of M-tame functions respectively
Landau-Ginzburg models. They are given precisely in
\cite[pages 589 and 590]{CV93}. They use $tt^*$-geometry.

Indeed, any matrix $S\in T(n,\R)$ gives together with 
arbitrary values $(u_1,...,u_n)$ with $u_i\neq u_j$
for $i\neq j$ and a sufficiently generic value $\xi\in S^1$
rise to a TERP-structure in the sense of \cite{He03},
more precisely, it gives a semisimple mixed TERP-structure
of weight 1 \cite[Lemma 10.1]{HS07},
which we call now $TERP(S,(u_1,...,u_n),\xi)$.

But for conjecture \ref{t1.9}, Cecotti and Vafa want to 
consider a limit TERP-structure for $(u_1,...,u_n)\to (0,...,0)$.
This should be the UV limit (ultraviolet limit). 
They {\it assume}
that it exists and that it has good properties, especially
it should be pure and polarized and have the correct spectrum.
In \cite[ch. 5, page 601]{CV93}, they {\it conclude} 
that the UV limit is well defined and nondegenerate 
(in a certain sense), 
if $S^{-1}S^t$ is semisimple.

We agree neither with the {\it assumption} nor with the 
{\it conclusion}.
The following example serves for both as a counterexample.

Therefore we do not consider conjecture \ref{t1.9}
(for the M-tame case respectively the Landau-Ginzburg models)
as proved in \cite{CV93}.
Though we do believe that $tt^*$-geometry is a promising
road. But a much more precise analysis of the limit 
behaviour seems to be needed.

\begin{example}\label{t8.1}
Consider a family of exceptional unimodal singularities, 
e.g. the family $E_{12}$: 
\begin{eqnarray}\label{8.1}
f_{t_\mu}(x,y)=x^3+y^7+t_\mu\cdot xy^5\quad\textup{with }
\mu=12.
\end{eqnarray}
$f_0$ is quasihomogeneous of weighted degree 1 with respect
to the weights $(w_x,w_y)=(\frac{1}{3},\frac{1}{7})$, 
and $f_{t_\mu}$ for $t_\mu\neq 0$
is semiquasihomogeneous.

The TERP-structures $TERP(f_{t_\mu})$ were studied in
\cite[8.3 (C)]{He03}: There is a bound $r_2\in\R_{>0}$
such that $TERP(f_{t_\mu})$ is not pure for $|t_\mu|=r_2$,
it is pure and polarized for $|t_\mu|<r_2$,
and it is pure, but not polarized for 
$|t_\mu|>r_2$. 
The spectral numbers (from Steenbrink's MHS) 
are called $\alpha_1,...,\alpha_\mu$
and satisfy here 
\begin{eqnarray}
\alpha_j+\alpha_{\mu+1-j}=0,\hspace*{3.3cm}\nonumber \\
\alpha_1=\frac{-11}{21}<\frac{-1}{2}<\alpha_2=\frac{-8}{21}
<...<\alpha_{\mu-1}=\frac{8}{21}<\frac{1}{2}<\frac{11}{21}
=\alpha_\mu.\label{8.2}
\end{eqnarray}
The eigenvalues of the supersymmetric index $\QQ$ are for
$|t_\mu|\neq r_2$
\begin{eqnarray}\label{8.3}
\alpha_2,...,\alpha_{\mu-1}\textup{ and }
\pm\Bigl(1-\frac{|t_\mu|^2}{r_2^2}\Bigr)^{-1}
\Bigl(\alpha_1-\frac{|t_\mu|^2}{r_2^2}(-1-\alpha_1)\Bigr).
\end{eqnarray}
The last two eigenvalues of $\QQ$ tend for
$|t_\mu|\to 0$ to $\pm\alpha_1=\mp\frac{11}{21}$
and for $|t_\mu|\to\infty$ to 
$\pm (-1-\alpha_1)=\mp\frac{10}{21}$.

Now consider a universal unfolding 
\begin{eqnarray}\label{8.4}
F_t(x,y)=f_{t_\mu}(x,y)+\sum_{j=1}^{\mu-1}t_jm_j,\quad
t\in M\subset\C^\mu,
\end{eqnarray}
for suitable monomials $m_j$ with weighted degree 
$\deg_w(m_j)<1$. 
Here $M\subset\C^\mu$ is an open set which contains 
$\C^{\mu-1}\times\{0\}\cup \{(0,...,0)\}\times\C$
and which is invariant under the flow of the Euler field
$E=\sum_{j=1}^\mu \deg_w(t_j)\cdot t_j\frac{\paa}{\paa t_j}$.

Choose $\xi\in S^1$ and choose for any $(u_1,...,u_\mu)\in\C^\mu$
with $\Ree\Bigl(\frac{u_i-u_j}{\xi}\Bigr)\neq 0$ for $i\neq j$
a special distinguished system of paths:
They shall go straight in the direction $\xi$ to 
$\paa\Delta_\eta$ and then turn on 
$\paa\Delta_\eta$ to $\xi\cdot\eta$. The set 
\begin{eqnarray}\label{8.5}
\{t\in M &|& \textup{the critical values }u_1,...,u_\mu 
\textup{ of }F_t\textup{ satisfy }\\
&&\Ree \Bigl(\frac{u_i-u_j}{\xi}\Bigr)\neq 0
\textup{ for }i\neq j\}\nonumber
\end{eqnarray}
consists of finitely many regions, the Stokes regions.
Each Stokes region gives one $G_{sign,\mu}$-orbit
of Stokes matrices $S$. For $t$ in one region
\begin{eqnarray}\label{8.6}
TERP(F_t) = TERP(S,(u_1,...,u_\mu),\xi),
\end{eqnarray}
and rescaling $(u_1,...,u_\mu)$ to 
$(r\cdot u_1,...,r\cdot u_\mu)$ with $r>0,r\to 0$,
corresponds to moving $t$ along $-\Ree E$.
There are now two severe problems.

\begin{list}{}{}
\item[(I)]
For $t\in M$ in one region
$TERP(F_t)$ tends to $TERP(f_0)$ only if $t_\mu=0$.
If $t_\mu\neq 0$ then for $r\to 0$ $TEPP(F_t)$ approximates
$TERP(f_{t_\mu})$ 
for larger and larger $t_\mu$, so it will become pure, 
but not polarized, and the eigenvalues of its supersymmetric
index $\QQ$ will tend to $\alpha_2,...,\alpha_{\mu-1},
\pm(-1-\alpha_1)$.
\item[(II)]
The $\textup{Br}_\mu\ltimes G_{sign,\mu}$-orbit of all Stokes
matrices is infinite \cite{Eb16}.
The $G_{sign,\mu}$-orbits of the Stokes matrices
from the finitely many Stokes regions in $M$
form only a finite subset. 
For $S$ not in this subset it is not at all clear how
$TERP(S,(u_1,...,u_\mu),\xi)$ will behave for $r\to 0$.
\end{list}

Both problems show that the assumption and the conclusion
about existence and good properties of the UV limit
are not justified in the generality in which they are
claimed in \cite{CV93}.
\end{example}

\subsection{(Harmonic) vector bundles}\label{s8.2}
We hope that the conjectures \ref{t1.6}, \ref{t1.7}
and \ref{t1.9} are true and will be proved in the future.
The special cases of the HOR-matrices made crucial use
of the formulas \eqref{4.20}
$(-1)^k\cdot S^{-1}S^t = R^{mat}_{(k)}(S)^n$ for $k\in\{1,2\}$. 
They are special cases of the formulas \eqref{4.12}
$(-1)^k\cdot S^{-1}S^t = R^{mat}_{(k1)}\circ ...\circ 
R^{mat}_{(kn)}$. Here the matrices $R^{mat}_{(kj)}$
are obtained by a certain twist from matrices for 
Picard-Lefschetz transformations, and they are
companion matrices (remark \ref{t4.3}). 
We hope that the formulas \eqref{4.12} will be useful
for an approach to the conjectures \ref{t1.6},
\ref{t1.7} and \ref{t1.9}.

Certainly, it will also be useful to consider the flat 
vector bundle on $\C-\{u_1,...,u_n\}$ of rank $n$
whose monodromy is given by these matrices
$R^{mat}_{(kj)}$ at $u_j$ (for $j\in\{1,...,n\}$) 
and by $(-1)^kS^{-t}S$ at $\infty$.
The vector bundle whose monodromy is given by the
Picard-Lefschetz transformations and $(-1)^kS^{-t}S$
is very familiar, it arises as homology bundle
of a suitable function with $A_1$-singularities only.
We hope that the local monodromies given by the 
companion matrices $R^{mat}_{(kj)}$ will become useful
beyond the special case of HOR-matrices.

In the special case of HOR-matrices, the flat bundle
decomposes because of \eqref{4.20} into flat subbundles,
for each eigenvalue $\kappa$ of $R^{mat}_{(k)}(S)$ one.
In the semisimple case, these are flat line bundles.
Then one can understand the $\uuuu\beta$ and
$Sp(S)$ in terms of natural holomorphic extensions
of these line bundles on $\C-\{u_1,...,u_n\}$
to $\P^1\C$. 

But how this observation might extend to the general
case of arbitrary matrices $S\in T(n,\R)$ is not clear to us.
Possibly work on harmonic bundles, tame or wild at
$\{u_1,...,u_n,\infty\}$, by Biquard, Boalch, Mochizuki
and Sabbah might be useful. And this might have 
connections to the TERP structures.

\subsection{Thom-Sebastiani formulas}\label{s8.3}
In the case of isolated hypersurface singularities (short: ihs), 
an important technique for obtaining new ihs is, to consider the sum 
$f(x_0,...,x_m)+g(x_{m+1},...,x_{m+n+1})$ of two ihs
$f$ and $g$ in different variables. 
This is discussed in \cite[I.2.7]{AGV88} and reviewed 
(in notations closer to this paper) in \cite{GH17}.
There is a canonical isomorphism
\begin{eqnarray}\label{8.7}
\Phi:Ml(f+g,1)&\stackrel{\cong}{\longrightarrow} 
&Ml(f,1)\otimes Ml(g,1),\\
\textup{with } M(f+g)&\cong & M(f)\otimes M(g) \label{8.8}\\
\textup{and } 
L^{hnor}(f+g)&\cong& L^{hnor}(f)\otimes L^{hnor}(g).\label{8.9}
\end{eqnarray}
If $\uuuu{\delta}=(\delta_1,...,\delta_{\mu(f)})$
and $\uuuu{\gamma}=(\gamma_1,...,\gamma_{\mu(g)})$ are
distinguished bases of $f$ and $g$  with Stokes matrices
$S(f)$ and $S(g)$, then 
$$\Phi^{-1}(\delta_1\otimes \gamma_1,...,
\delta_1\otimes \gamma_{\mu(g)},
\delta_2\otimes \gamma_1,...,
\delta_2\otimes \gamma_{\mu(g)},
...,
\delta_{\mu(f)}\otimes \gamma_1,...,
\delta_{\mu(f)}\otimes \gamma_{\mu(g)})$$
is a distinguished basis of $Ml(f+g,1)$,
that means, one takes the vanishing cycles 
$\Phi^{-1}(\delta_i\otimes \gamma_j)$ in the lexicographic order.
Then by \eqref{6.7} and \eqref{8.9}, the matrix 
\begin{eqnarray}\label{8.10}
S(f+g)=S(f)\otimes S(g)
\end{eqnarray}
(where the tensor product is defined
so that it fits to the lexicographic order) 
is the Stokes matrix of this distinguished basis.

In \cite[ch. 8]{SS85} a Thom-Sebastiani for
Steenbrink's mixed Hodge structure is stated. 
It is fine if the monodromy is semisimple, but it 
needs a correction in the general case.
That correction is an interesting and nontrivial twist
\cite[Corollary 6.5]{BH17}, which comes from a
Fourier-Laplace transformation.
Anyway, the resulting Thom-Sebastiani formula in \cite{SS85}
for the spectral pairs of $f$, $g$ and $f+g$ is correct.

The set of HOR-matrices is not invariant under the
tensor product of matrices. It might be a good idea
to check whether there are natural modifications
for the recipe how the HOR-matrices give rise to spectral numbers,
which are compatible with the Thom-Sebastiani formulas.

\end{document}